\documentclass[11pt]{amsart}
\usepackage{amssymb,latexsym,amsmath,graphicx,graphics,epic,eepic}

\theoremstyle{plain}
\newtheorem{theorem}{Theorem}
\numberwithin{theorem}{section}
\newtheorem{lemma}[theorem]{Lemma}
\newtheorem{proposition}[theorem]{Proposition}
\newtheorem{corollary}[theorem]{Corollary}

\theoremstyle{definition}

\newtheorem{remark}[theorem]{Remark}

\newcommand{\C}{{\mathbb C}}
\newcommand{\D}{{\mathbb D}}
\newcommand{\R}{{\mathbb R}}
\newcommand{\Z}{{\mathbb Z}}
\newcommand{\Q}{{\mathbb Q}}
\renewcommand{\P}{{\mathbb P}}
\newcommand{\s}{{\mathbb S}}

\begin{document}
\title{Symmetries and exotic smooth structures on a $K3$ surface}
\author{Weimin Chen and Slawomir Kwasik}
\thanks{The first author is supported in part by NSF grant DMS-0603932}
\date{}

\begin{abstract}
Smooth and symplectic symmetries of an infinite family of distinct 
exotic $K3$ surfaces are studied, and comparison with the corresponding 
symmetries of the standard $K3$ is made. The action on the $K3$ lattice 
induced by a smooth finite group action is shown to be strongly restricted, 
and as a result, nonsmoothability of actions induced by a holomorphic 
automorphism of a prime order $\geq 7$ is proved and nonexistence 
of smooth actions by several $K3$ groups is established (included among
which is the binary tetrahedral group $T_{24}$ which has the smallest order). 
Concerning symplectic symmetries, the fixed-point set structure of a 
symplectic cyclic action of a prime order $\geq 5$ is explicitly determined, 
provided that the action is homologically nontrivial.
\end{abstract}

\maketitle

\section{Introduction}

The main purpose of this paper is to investigate the effect of a change 
of a smooth structure on the smooth symmetries of a closed, oriented 
$4$-dimensional smoothable manifold. The influence of symmetries on smooth
structures on a manifold is one of the basic questions in the theory 
of differentiable transformation groups. The following classical theorem
of differential geometry gives a beautiful characterization of the standard 
sphere $\s^n$ among all simply connected manifolds. It led to an extensive
study of various degrees of symmetry for the (higher dimensional)
exotic spheres in the 1960s and
70s (cf. \cite{Hsiang}). Lawson and Yau even found that there exist
exotic spheres which support no actions of small groups such as $\s^3$ or
$SO(3)$ (cf. \cite{LY}). See \cite{Sch} for a survey. 
  
\vspace{3mm}

{\bf Theorem} (A Characterization of $\s^n$). 
{\it Let $M^n$ be a closed, simply connected 
manifold of dimension $n$, and let $G$ be a compact Lie group which
acts smoothly and effectively on $M^n$. Then $\dim G\leq n(n+1)/2$,
with equality if and only if $M^n$ is diffeomorphic to $\s^n$.}

\vspace{3mm}

The subject of symmetries of exotic smooth $4$-manifolds, on the other hand,
has been so far rather an untested territory.
Our investigations of smooth symmetries of $4$-manifolds have been 
focused on the case of $K3$ surfaces. These manifolds exhibit surprisingly 
rich geometric structures and have been playing one of the central 
roles in both the theory of complex surfaces and topology of smooth
$4$-manifolds. 

To be more specific, we will study symmetries of an infinite family 
of distinct, closed, oriented smooth 4-manifolds, each of which is 
orientation-preservingly homeomorphic, but not diffeomorphic, to a $K3$ 
surface (canonically oriented as a complex surface). These exotic $K3$ 
surfaces, originally due to Fintushel and Stern, are obtained by performing 
the knot surgery construction simultaneously on three disjoint 
embedded tori representing distinct homology classes in a Kummer surface 
(cf. \cite{FS}, compare also \cite{GM}). It is known that none of 
these $4$-manifolds support a complex structure (cf. \cite{FS, GS}), however, 
one may arrange the knot surgeries so that each of these manifolds supports 
a symplectic structure compatible with the given orientation (cf. \cite{FS}). 

A $K3$ surface is a simply-connected complex surface with trivial 
canonical bundle. It is known that all $K3$ surfaces are deformation
equivalent as complex surfaces (therefore diffeomorphic as oriented smooth
$4$-manifolds), and that all $K3$ surfaces are K\"{a}hler surfaces 
(cf. \cite{BPV}). There is an extensive study on finite subgroups of the
automorphism group of a $K3$ surface, beginning with the fundamental work 
of Nikulin \cite{N}. Special attention has been given to those subgroups
of automorphisms which induce a trivial action on the canonical bundle of
the $K3$ surface. (Such automorphisms are called symplectic; in Nikulin 
\cite{N} they were called algebraic.) A finite group $G$ is called a $K3$ 
group (resp. symplectic $K3$ group) if $G$ can be realized as a subgroup of 
the automorphism group (resp. symplectic automorphism group) of a  $K3$ 
surface. Finite abelian groups of symplectic automorphisms of a $K3$ surface
were first classified by Nikulin in \cite{N}; in particular it was shown 
that a finite symplectic automorphism must have order $\leq 8$. Subsequently,
Mukai \cite{Mu} determined all the symplectic $K3$ groups (see also 
\cite{Kon, Xiao}). There are $11$ maximal ones, all of which are 
characterized as certain subgroups of the Mathieu group $M_{23}$. Finally,
a cyclic group of prime order $p>7$ is a $K3$ group (necessarily 
non-symplectic) if and only if $p\leq 19$ (cf. \cite{N, MO}).

We recall three relevant properties of automorphism groups of $K3$ 
surfaces. First, a finite-order automorphism of a $K3$ surface preserves 
a K\"{a}hler structure, hence by the Hodge theory, it is symplectic if and 
only if the second cohomology contains a $3$-dimensional subspace consisting
of invariant elements of positive square. Secondly, since 
a symplectic automorphism acts trivially on the canonical bundle, it follows 
that the induced representation at a fixed point (called a local 
representation) lies in $SL_2(\C)$; in particular, the fixed point is 
isolated. (Such actions are called pseudofree.) Finally, a nontrivial 
automorphism of a $K3$ surface must act nontrivially on the homology (cf. 
\cite{BPV}). 

Finite groups of automorphisms of a $K3$ surface are primary sources of
smooth and symplectic symmetries of the standard $K3$. 
(In fact, no examples of smooth symmetries of the standard $K3$ 
are known to exist that are not automorphisms of a $K3$ surface.)
Thus in analyzing symmetry properties of an exotic $K3$ surface, we will 
use these automorphisms as the base of our comparison. 

We shall now state our main theorems. In what
follows, we will denote by $X_\alpha$ a member of the infinite family of 
exotic $K3$ surfaces of Fintushel and Stern we alluded to earlier. 
(A detailed review of their construction along with some relevant properties 
will be given in Section $2$; we point out here that the index 
$\alpha$ stands for a triple $(d_1,d_2,d_3)$ of integers which obey 
$1<d_1<d_2<d_3$ and are pairwise relatively prime.) 

The induced action on the quadratic form and the fixed-point set structure 
are two fundamental pieces of information associated with a finite group
action on a simply-connected $4$-manifold. In this regard, we have  

\begin{theorem}
Let $G\equiv \Z_p$ where $p$ is an odd prime. The following statements are
true for a smooth $G$-action on $X_\alpha$.
\begin{itemize}
\item [{(1)}] The induced action is trivial on a $3$-dimensional subspace of
$H^{2}(X_\alpha;\R)$ over which the cup-product is positive definite. 
\item [{(2)}] For $p\geq 7$, there is a $G$-invariant, orthogonal 
              decomposition of the intersection form on $H_2(X_\alpha;\Z)$ 
              as 
$$
3\left(\begin{array}{cc}
0 & 1\\
1 & 0\\
\end{array} \right)
\oplus 2(-E_8)
$$  
where the induced $G$-action on each hyperbolic summand is trivial.
\end{itemize}
\end{theorem}

\begin{remark}
(1) For a smooth $\Z_p$-action on a homotopy $K3$ surface,
Theorem 1.1 holds true automatically when $p\geq 23$ because in this case 
the action is necessarily homologically trivial. However, when $p<23$,
nothing is known in general about the induced action on the $K3$ lattice.
For a prime order symplectic automorphism of a $K3$ surface, Nikulin showed
in \cite{N} that the action on the $K3$ lattice is unique up to conjugacy,
which was explicitly determined in \cite{Mo, GarSa} by examing some concrete
examples of the action. In particular, Theorem 1.1 (2) is false for symplectic 
automorphisms of a $K3$ surface (cf. the proof of Corollary 1.3).

(2) The $G$-invariant, orthogonal decomposition in Theorem 1.1 (2) gives
severe restrictions on the induced integral $G$-representation on 
$H_2(X_\alpha;\Z)$; in particular, when $p>7$ the action must be 
homologically trivial because $\text{Aut }(E_8)$ contains no elements of
order $>7$. Note that one does not expect such a result in general, as for
each prime $p$ with $p<23$, there exists an automorphism of a $K3$ surface
of order $p$, which is necessarily homologically nontrivial.

(3) Let $F$ be the fixed-point set.  A general result of Edmonds 
(cf. Proposition 3.1) constrains the number of $2$-dimensional components of 
nonzero genus in $F$ by equating the first Betti number of $F$ (if nonempty) and
the number of summands of cyclotomic type in the induced
integral $G$-representation on the second homology. For a smooth 
$\Z_p$-action on a homotopy $K3$ in general, there are no summands of 
cyclotomic type when $p\geq 13$, and consequently $F$ does not contain any
$2$-dimensional non-spherical components in these cases. However, when $p=7$ 
or $11$, such a summand does occur. In fact, for both $p=7$ and $p=11$, 
there exists an automorphism of a $K3$ surface of order $p$ which fixes a 
regular fiber of an elliptic fibration on the $K3$ surface (cf. \cite{MO}). 
With the above observations, note that Theorem 1.1 (2) implies that for a 
smooth $\Z_p$-action on $X_\alpha$ of order $p=7$ or $11$, $F$ contains 
at most $2$-dimensional spherical components (cf. Lemma 4.5), which is in 
contrast with the case of the standard $K3$ we mentioned earlier. Finally, 
a calculation with the Lefschetz fixed point theorem indicates that for 
$p\geq 7$, $F$ also has a fairly large size, e.g., $\chi(F)\geq 10$. 
(In contrast a symplectic automorphism of a $K3$ surface of order $7$ has 
only three isolated fixed points, hence $\chi(F)=3$.) 

\end{remark}

We have seen from the discussions in Remark 1.2 that for any prime $p\geq 7$, 
a smooth $\Z_p$-action on $X_\alpha$ differs in many aspects from an 
automorphism of a $K3$ surface. In particular, we note the following 
relative nonsmoothability result as a corollary of Theorem 1.1.

Recall that each $X_\alpha$ is homeomorphic to a $K3$ surface. Thus any
finite-order automorphism of a $K3$ surface induces a locally linear 
topological action on $X_\alpha$ after we fix a homeomorphism between 
$X_\alpha$ and the standard $K3$.

\begin{corollary}
Any locally linear topological action induced by an automorphism 
of a $K3$ surface of a prime order $\geq 7$ is nonsmoothable on 
$X_\alpha$.
\end{corollary}

\begin{proof}
Let $g$ be an automorphism of a $K3$ surface of a prime order $p\geq 7$.
If $g$ is non-symplectic, then $g$ is not smoothable on $X_\alpha$ 
by Theorem 1.1 (1). Suppose $g$ is a symplectic automorphism. Then 
by Nikulin \cite{N}, we have $p=7$, and moreover, the action of 
$g$ is pseudofree with $3$ isolated fixed points. Suppose $g$ is smoothable
on $X_\alpha$. Then by Theorem 1.1 (2)
and Lemma 4.5, the trace of the action of $g$ on $H_2(X_\alpha;\Z)$ is
at least $8$, so that by the Lefschetz fixed point theorem (cf. Theorem 3.4),
the Euler number of the fixed point set of $g$ is at least $10$.
A contradiction. 

\end{proof}

Next we turn our attention to smooth involutions, i.e., 
smooth $\Z_2$-actions on $X_\alpha$. Let $g:X_\alpha
\rightarrow X_\alpha$ be any smooth involution. Since $X_\alpha$ is 
simply-connected, $g$ can be lifted to the spin bundle over $X_\alpha$,
where there are two cases: (1) $g$ is of even type, meaning that the 
order of lifting to the spin bundle is $2$, or (2) $g$ is of odd type, 
meaning that the order of lifting to the spin bundle is $4$. 
Moreover, $g$ has $8$ isolated fixed points in the case of an even type, 
and $g$ is free or has only $2$-dimensional fixed components in the case 
of an odd type (cf. \cite{AB, Bry}).

\begin{theorem}
Suppose $g:X_\alpha\rightarrow X_\alpha$ is an odd type smooth 
involution. Then the fixed-point set of $g$ belongs to one of the following 
three possibilities:
\begin{itemize}
\item [{(1)}] An empty set.
\item [{(2)}] A disjoint union of two tori.
\item [{(3)}] A disjoint union of spheres or tori where the number of tori
is at most one.
\end{itemize}
\end{theorem}

Let $\tau$ be an anti-holomorphic involution on a $K3$ surface. (Note that
$\tau$ is holomorphic with respect to some other complex structure on the
smooth $4$-manifold, cf. \cite{D}.) Then $\tau$ falls into one of the 
following three types according to the fixed point set $\text{Fix}(\tau)$ 
of $\tau$ (cf. \cite{N1}); in particular, $\tau$ is of odd type:
\begin{itemize}
\item $\text{Fix}(\tau)=\emptyset$, 
\item $\text{Fix}(\tau)$ is a union of two tori, 
\item $\text{Fix}(\tau)$ is a union of orientable surfaces of genus $\leq 10$,
such that the number of non-spherical components in $\text{Fix}(\tau)$
is at most one.
\end{itemize}

\begin{corollary}
A locally linear topological $\Z_2$-action induced by an anti-holomorphic 
involution is nonsmoothable on $X_\alpha$ if it has 
a fixed component of genus $\geq 2$. {\em (}Note that such 
anti-holomorphic involutions do exist, cf. \cite{N1}.{\em)}
\end{corollary}

\begin{remark}
There are previously known examples of locally linear topological actions
on closed $4$-manifolds which are not smoothable. For example, there is a 
locally linear, pseudofree, homologically trivial topological action of 
order $5$ on $\C\P^2\#\C\P^2$ which can not be realized as an equivariant 
connected sum of two copies of $\C\P^2$ (cf. \cite{EE}). By the main result 
of \cite{HL}, the action is not smoothable with respect to any smooth 
structure on $\C\P^2\#\C\P^2$. (For more recent examples of nonsmoothable 
actions on closed $4$-manifolds, see e.g. \cite{LN1, LN2}, 
and for nonsmoothable 
actions on non-closed $4$-manifolds, see \cite{KL}.) However, the 
nonsmoothability in Corollary 1.3 and Corollary 1.5 is of a different 
nature; the action is smooth (even holomorphic) for one smooth structure 
but not smoothable with respect to some (in fact infinitely many) other 
smooth structures. 

\end{remark}

Our investigation of the possible effect of a change of smooth structures 
on the smooth symmetries of a closed, oriented smoothable $4$-manifold 
is based on the following simple fact. Suppose $M^4$ is a simply-connected,
oriented smooth $4$-manifold with an orientation-preserving smooth action 
of a finite group $G$. Let $L$ be the primitive sublattice of $H^2(M^4;\Z)$ 
spanned by the Seiberg-Witten basic classes of $M^4$ (we assume 
$b_2^{+}(M^4)>1$). Then the induced $G$-action on $H^2(M^4;\Z)$ preserves 
$L$ as it preserves the set of Seiberg-Witten basic classes. One can try 
to analyze the $G$-action on $H^2(M^4;\Z)$ through the actions on $L$ and
$L^\perp$, the orthogonal complement of $L$ in $H^2(M^4;\Z)$. With this
understood, a crucial ingredient in our investigation is the following
property of $X_\alpha$: $L$ is isotropic and of rank $3$, such that 
$$
L^\perp/L=2(-E_8). 
$$ 
(See Lemma 4.2 for more details.) Furthermore, one can arrange $X_\alpha$
such that each Seiberg-Witten basic class is fixed under the action up to
a change of signs; in particular, an odd order $G$ must act trivially on $L$. 
Given this, Theorem 1.1 follows readily by analyzing the corresponding 
action on $L^\perp/L=2(-E_8)$ where $G$ is cyclic of a prime order $\geq 7$. 

The above mentioned property of $X_\alpha$ 
can be further exploited to prove non-existence of effective
smooth $G$-actions on $X_\alpha$ for a certain kind of finite groups $G$. For
instance, suppose $G$ is of odd order and there are no nontrivial $G$-actions
on the $E_8$ lattice (e.g. $G$ is a $p$-group with $p>7$), 
then any smooth $G$-actions on $X_\alpha$ must be 
homologically trivial, and therefore, by a theorem of McCooey \cite{McC}
$G$ must be abelian of rank at most $2$ (cf. Corollary 4.4). In particular,
a nonabelian $p$-group with $p>7$ can not act smoothly and effectively on
$X_\alpha$. 
We remark that while for a given finite group $G$, we do not know a priori any 
obstructions to the existence of a smooth $G$-action on a homotopy $K3$
surface, 
a non-existence result of this sort for $X_\alpha$ may be in fact 
purely topological in nature. 

The following non-existence theorem of smooth actions on $X_\alpha$ 
covers the cases of several $K3$ groups, hence it must not be purely 
topological in nature. Its proof requires a deeper analysis
of the induced actions on the $E_8$ lattice. 

\begin{theorem}
Let $G$ be a finite group whose commutator $[G,G]$ contains a subgroup
isomorphic to $(\Z_2)^4$ or $Q_8$, where in the case of $Q_8$
the elements of order $4$ in the subgroup are conjugate in $G$. 
Then there are no effective smooth $G$-actions on $X_\alpha$. 
\end{theorem}

A complete list of symplectic $K3$ groups along with their commutators
can be found in Xiao \cite{Xiao}, Table $2$. By examing the list we note
that among the $11$ maximal $K3$ groups the following can not act smoothly
and effectively on $X_\alpha$ (cf. Corollary 5.4): 
$$
M_{20}, F_{384}, A_{4,4}, T_{192}, H_{192}, T_{48}.
$$
We also note that the binary tetrahedral group $T_{24}$ of order $24$ is 
the $K3$ group of the smallest order which can not act smoothly and 
effectively on $X_\alpha$ by Theorem 1.7. 

\vspace{3mm}

As we mentioned earlier, each exotic $K3$ $X_\alpha$ supports an 
orientation-compatible symplectic structure. In order to investigate 
how symplectic symmetries may depend on the underlying smooth
structure of a $4$-manifold, we also analyzed finite group actions 
on $X_\alpha$ which preserve an orientation-compatible symplectic 
structure.

Recall that $\text{Aut}(E_8\oplus E_8)$ is a semi-direct
product of $\text{Aut}(E_8)\times \text{Aut}(E_8)$ by $\Z_2$ (cf. \cite{Se}).
Thus for any smooth $G$-action on $X_\alpha$, where $G\equiv\Z_p$ for an
odd prime $p$, there is an associated homomorphism 
$$
\Theta=(\Theta_1,\Theta_2):G=\Z_p\rightarrow \text{Aut}(E_8)\times 
\text{Aut}(E_8).
$$ 
The following theorem gives a complete description of the fixed-point
set structure of a symplectic $\Z_p$-action on $X_\alpha$ for $p\geq 5$,
provided that both homomorphisms $\Theta_1,\Theta_2:\Z_p\rightarrow
\text{Aut}(E_8)$ are nontrivial. Note that this implies that $p=5$ or 
$p=7$.

\begin{theorem}
Let $F\subset X_\alpha$ be the fixed-point set of a symplectic 
$\Z_p$-action of a prime 
order $p\geq 5$ such that both $\Theta_1,\Theta_2$ are nontrivial. 
Set $\mu_p\equiv \exp(\frac{2\pi i}{p})$. Then 

(1) if $p=5$, there are two possibilities:
\begin{itemize}
\item [{(i)}] $F$ consists of $14$ isolated fixed points, two with local 
representation $(z_1,z_2)\mapsto (\mu_p^k z_1,\mu_p^k z_2)$, six with
local representation $(z_1,z_2)\mapsto (\mu_p^{4k} z_1,\mu_p^{3k} z_2)$,
two with local representation 
$(z_1,z_2)\mapsto (\mu_p^{2k} z_1,\mu_p^{4k} z_2)$, and four with
local representation $(z_1,z_2)\mapsto (\mu_p^{3k} z_1,\mu_p^{3k} z_2)$,
for some $k\neq 0\pmod{p}$. 
\item [{(ii)}] $F=F_1\sqcup F_2$, where $F_1$ consists of $4$ isolated 
fixed points with local representations 
$$
(z_1,z_2)\mapsto (\mu_p^q z_1,\mu_p^{-q} z_2), \mbox{ evaluated at } 
q=1,2,3,4,
$$
and $F_2$ is divided into groups of fixed points of the following two types
where the number of groups is less than or equal to $2$ {\em(}in particular, 
$F_2$ may be empty{\em)}:
\begin{itemize}
\item $3$ isolated points, two with local representation 
$(z_1,z_2)\mapsto (\mu_p^{-3k} z_1,\mu_p^{-k} z_2)$ and one with local
representation $(z_1,z_2)\mapsto (\mu_p^{3k} z_1,\mu_p^{3k} z_2)$, and
one fixed $(-2)$-sphere with local representation $z\mapsto \mu_p^k z$,
for some $k\neq 0\pmod{p}$,
\item a fixed torus of self-intersection $0$.
\end{itemize}
\end{itemize}

(2) if $p=7$, $F$ consists of $10$ isolated fixed points, two with local
representation $(z_1,z_2)\mapsto (\mu_p^{2k} z_1,\mu_p^{3k} z_2)$, two with
local representation $(z_1,z_2)\mapsto (\mu_p^{-k} z_1,\mu_p^{-k} z_2)$,
two with local representation 
$(z_1,z_2)\mapsto (\mu_p^{2k} z_1,\mu_p^{4k} z_2)$, and four with local
representation $(z_1,z_2)\mapsto (\mu_p^{-2k} z_1,\mu_p^{k} z_2)$,
for some $k\neq 0\pmod{p}$.
\end{theorem}

\begin{remark}

(1) We remark that by the work of Edmonds and Ewing \cite{EE}, the fixed-point
set structure of a pseudofree action in Theorem 1.8 can be actually 
realized by a locally linear, topological action on $X_\alpha$. On the other
hand, none of the known obstructions to smoothability of topological actions 
(see Section $3$) could rule out the possibility that the fixed-point set structure 
may be realized by a smooth or even symplectic action. 

(2) Since the case of small primes $p$ is missing in Theorem 1.1, Theorem
1.8 can be viewed as a complement to Theorem 1.1. We remark that the 
homological triviality of actions in Theorem 1.1 for the case of $p>7$ 
plus the detailed information about the (homologically nontrivial) $\Z_5$ 
and $\Z_7$ actions in Theorem 1.8 put considerable limitations on the 
symplectic actions of an arbitrary finite group on $X_\alpha$. 

\end{remark}

The proof of Theorem 1.8 is based on a combination of the pseudoholomorphic 
curve techniques developed in our previous work \cite{CK} and a delicate
exploitation of the induced actions on the $E_8$ lattice. Note that the 
latter is possible only because of the property $L^\perp/L=2(-E_8)$ of 
$X_\alpha$. For a general homotopy $K3$ surface, a symplectic $\Z_p$-action 
could have a much more complicated fixed-point set structure. 
However, if the finite 
group which acts symplectically on a homotopy $K3$ has a relatively 
complicated group structure (e.g., a maximal symplectic $K3$ group),
then the fixed-point set structure can also be explicitly determined.
This observation was systematically exploited in our subsequent paper 
\cite{CK1} where the following problem was investigated. 

\vspace{3mm}

{\bf Problem} \hspace{2mm} {\it Let $X$ be a homotopy $K3$ surface supporting
an effective action of a ``large'' $K3$ group via symplectic symmetries.
What can be said about the smooth structure on $X$?}

\vspace{3mm}

In particular, a characterization of the ``standard'' smooth structure of
$K3$ in terms of symplectic symmetry groups was obtained (compare with the 
corresponding characterization of $\s^n$ at the beginning of the 
introduction). See \cite{CK1} for more details.

\vspace{3mm}

The current paper is organized as follows.

In Section $2$ we give a detailed description of the Fintushel-Stern 
exotic $K3$'s that are to be considered in this paper, along with their 
relevant properties.

In Section $3$ we collect various known results concerning topological and
smooth actions of finite groups on $4$-manifolds. These results are used
in our paper (sometimes successfully and sometimes not) to measure the 
difference between the symmetries of the standard and exotic $K3$ surfaces.
In particular, these results are the criteria used in the proof of 
Theorem 1.8, with which the fixed-point set structure of the group action is 
analyzed. 

Sections $4$, $5$ and $6$ contain proofs of Theorem 1.1, Theorem 1.4,
Theorem 1.7 and Theorem 1.8.

\section{The Fintushel-Stern exotic $K3$'s}

The construction of this type of exotic $K3$'s was briefly mentioned in the 
paper of Fintushel and Stern \cite{FS}. In this section we give a detailed
account of one particular family of such exotic $K3$'s that are used in this
paper, along with proofs of some relevant properties that will be used in
later sections.

The exotic $K3$ surfaces are the $4$-manifolds that result from performing
the knot surgery construction in \cite{FS} simultaneously on three disjoint 
embedded tori in a Kummer surface. We begin with a topological
description of a Kummer surface (following \cite{GM}) and establish 
some relevant properties of the three disjoint tori in it. 

Let $\s^1$ be the unit circle in $\C$. Let $T^4$ denote the $4$-torus
$\s^1\times\s^1\times\s^1\times\s^1$ and $\rho:\s^1\rightarrow \s^1$ denote 
the complex conjugation respectively. Moreover, we shall let $\hat{\rho}$
denote the corresponding diagonal involution on $T^4$ or $T^2\equiv \s^1\times \s^1$.
Then the underlying $4$-manifold $X$ 
of a Kummer surface is obtained by replacing each of the $16$ 
singularities $(\pm 1,\pm 1,\pm 1, \pm 1)$ 
in $T^4/\langle \hat{\rho}\rangle$ by a $(-2)$-sphere.
More precisely, for each of the $16$ singularities we shall remove a 
regular neighborhood of it and then glue back a regular neighborhood of
an embedded $(-2)$-sphere (which abstractly is a $D^2$-bundle over $\s^2$
with Euler number $-2$). Since the gluing is along $\R\P^3$ which has the
property that a self-diffeomorphism is isotopic to identity if and only if
it is orientation-preserving (cf. \cite{B}), the resulting $4$-manifolds 
for different choices of the gluing map are diffeomorphic to each other. 
In fact, they can be identified by a diffeomorphism which is identity on 
$T^4/\langle \hat{\rho}\rangle$ with a regular neighborhood of the
$16$ singularities removed and sends the corresponding embedded $(-2)$-spheres
diffeomorphically onto each other. Our $4$-manifold $X$ is simply a fixed
choice of one of these $4$-manifolds. As for the orientation of $X$, we shall
orient $T^4$ by $d\theta_0\wedge d\theta_1\wedge d\theta_2\wedge d\theta_3$,
where $\theta_j$, $j=0,1,2,3$, is the angular coordinate (i.e. 
$z=\exp(i\theta)$, $z\in\s^1$) on the $(j+1)$-th copy of $\s^1$ in $T^4$,
and the manifold $X$ is oriented by the orientation on 
$T^4/\langle \hat{\rho}\rangle$, whose smooth part is contained in
$X$.  

For $j=1,2,3$, let $\pi_j:T^4/\langle \hat{\rho}\rangle\rightarrow
T^2/\langle \hat{\rho}\rangle$ be the map induced by the projection
$$
(z_0,z_1,z_2,z_3)\mapsto (z_0,z_j).
$$
There is a complex structure $J_j$ on $T^4$, which is 
compatible with the given orientation on $T^4$, such that 
$\pi_j:T^4/\langle \hat{\rho}\rangle\rightarrow
T^2/\langle \hat{\rho}\rangle$ is holomorphic. 
Let $X(j)$ be the minimal complex surface obtained by 
resolving the singularities of $T^4/\langle \hat{\rho}\rangle$.
Then $\pi_j$ induces an elliptic fibration $X(j)\rightarrow T^2/\langle \hat{\rho}\rangle\equiv\s^2$. After
fixing an identification between $X(j)$ and $X$ in the manner described in
the preceding paragraph, we obtain three $C^\infty$-elliptic fibrations 
(cf. \cite{FM}) $\pi_j: X\rightarrow\s^2$. 

Given this, the three disjoint tori in $X$ which will be 
used in the knot surgery are some fixed regular fibers $T_j=\pi^{-1}_j
(\delta_j,i)$ of $\pi_j: X\rightarrow\s^2$, where $\delta_j\in \s^1$, 
$j=1,2,3$, are not $\pm 1$ and are chosen so that their images are distinct 
in $\s^1/\langle\rho\rangle$. (Note that $T_1$, $T_2$, $T_3$ are disjoint because the 
$z_0$-coordinates $\delta_1,\delta_2,\delta_3$ have distinct images in 
$\s^1/\langle\rho\rangle$.)

Concerning the relevant properties of the tori $T_1$, $T_2$ and $T_3$, we
first observe 

\begin{lemma}
The three disjoint tori $T_1$, $T_2$ and $T_3$ have the following
properties.
\begin{itemize}
\item [{(1)}] There are homology classes $v_1,v_2,v_3\in H_2(X;\Z)$
such that $v_i\cdot [T_j]=1$ for $i=j$ and $v_i\cdot [T_j]=0$ otherwise.
In particular, $[T_1]$, $[T_2]$, $[T_3]$ are all primitive classes and span
a primitive sublattice of rank $3$ in $H_2(X;\Z)$. 
\item [{(2)}] There are orientation compatible symplectic structures on $X$
with respect to which $T_1$, $T_2$ and $T_3$ are symplectic submanifolds.
\end{itemize}
\end{lemma}

\begin{proof} Observe that for each torus $T_j$, there is a sphere $S_j$
in the complement of the other two tori in $X$ which intersects $T_j$
transversely at a single point. For instance, for the torus $T_1$,
$S_1$ may be taken to be the proper transform of the section 
$$T^2\times\{1\}\times\{1\}/\langle\hat{\rho}\rangle$$ of 
the fibration $\pi_1:T^4/\langle \hat{\rho}\rangle\rightarrow
T^2/\langle \hat{\rho}\rangle$ 
in the complex surface $X(1)$. 
Here $S_1$ is regarded as a sphere in $X$ 
under the fixed identification between $X(1)$ and $X$. Part (1) of the 
lemma follows immediately.

Next we show that there are orientation compatible symplectic structures 
on $X$ with respect to which all three tori $T_1, T_2$ and $T_3$ are
symplectic. To see this, let $\theta_j$, $j=0,1,2,3$, be the
angular coordinate (i.e. $z=\exp(i\theta)$, $z\in\s^1$) on the $(j+1)$-th
copy of $\s^1$ in $T^4$. Then the following is a symplectic
$2$-form on $T^4$ which is equivariant with respect to the
diagonal involution $\hat{\rho}$:
$$
\sum_{(i,j,k)}(d\theta_0\wedge d\theta_i+d\theta_j\wedge
d\theta_k)
$$
where the sum is over $(i,j,k)=(1,2,3), (2,3,1)$ and $(3,1,2)$.
This gives rise to a symplectic structure on the orbifold
$T^4/\langle\hat{\rho}\rangle$. One can further
symplectically resolve the orbifold singularities to obtain a
symplectic structure on $X$ as follows. By the equivariant
Darboux' theorem, the symplectic structure is standard near each
orbifold singularity. In particular, it is modeled on a neighborhood of 
the origin in $\C^2/\{\pm 1\}$ and admits a Hamiltonian $\s^1$-action with 
moment map $\mu:(w_1,w_2)\mapsto\frac{1}{4}(|w_1|^2+|w_2|)$, where $w_1,w_2$
are the standard coordinates on $\C^2$. Fix a sufficiently small $r>0$
and remove $\mu^{-1}([0,r))$ from $T^4/\langle\hat{\rho}\rangle$
at each of its singular point. Then $X$ is diffeomorphic to the
$4$-manifold obtained by collapsing each orbit of the Hamiltonian
$\s^1$-action on the boundaries $\mu^{-1}(r)$, which is naturally a
symplectic $4$-manifold (cf. \cite{Le}). It is clear from the construction 
that all three tori $T_1, T_2$ and $T_3$ are symplectic, and moreover, the
symplectic structures are orientation compatible. 

\end{proof}

Following \cite{FS}, we call any Laurent polynomial
$$
P(t)=a_0+\sum_{j=1}^n a_j(t^j+t^{-j})
$$
in one variable with coefficient sum
$$
a_0+2\sum_{j=1}^n a_j=\pm 1, a_j\in\Z
$$
an $A$-polynomial. According to \cite{FS}, given any three $A$-polynomials 
$P_1(t)$, $P_2(t)$, $P_3(t)$, one can perform the so-called knot surgeries 
simultaneously along the tori $T_1, T_2, T_3$ to obtain an oriented 
$4$-manifold $X_{P_1P_2P_3}$, which is orientation-preservingly homeomorphic 
to $X$ and has Seiberg-Witten invariant
$$
SW_{X_{P_1P_2P_3}}=P_1(t_1)P_2(t_2)P_3(t_3),
$$
where $t_j=\exp(2[T_j])$, $j=1,2,3$. We remark that the homology
classes of $X_{P_1P_2P_3}$ are naturally identified with those
of $X$, and here $[T_j]$ in $t_j=\exp(2[T_j])$ denotes the
Poincar\'{e} dual of the class in $H_2(X_{P_1P_2P_3};\Z)$
which corresponds to the class of the torus $T_j$ in $X$ under the
identification. (In this paper, we shall use $[T_j]$ to denote either 
the homology class of the torus $T_j$ or the cohomology class that is 
Poincar\'{e} dual to $T_j$. The actual meaning is always clear from the 
context.) Moreover, when $P_1(t)$, $P_2(t)$, $P_3(t)$ are monic (i.e.,
the coefficient $a_n=\pm 1$), the $4$-manifold $X_{P_1P_2P_3}$ admits
orientation compatible symplectic structures because of Lemma 2.1 (2). 

We shall consider one particular infinite family of $(P_1(t),P_2(t),P_3(t))$
where each $A$-polynomial is monic and has the form 
$$
P_j(t)=1-(t^{d_j}+t^{-d_j}), j=1,2,3.
$$
Here $d_1,d_2,d_3$ are integers which obey $1<d_1<d_2<d_3$ and are pairwise 
relatively prime. We denote the corresponding $4$-manifold $X_{P_1P_2P_3}$ by 
$X(d_1,d_2,d_3)$. 

\begin{lemma}
For any orientation compatible symplectic structure $\omega$ on 
$X(d_1,d_2,d_3)$, one has $[T_j]\cdot [\omega]\neq 0$ for all $j$. 
If we assume {\em(}without loss of generality{\em)} that 
$[T_j]\cdot [\omega]> 0$ for all $j$, then the canonical class is given by
$$
c_1(K)=2\sum_{j=1}^3 d_j[T_j]. 
$$
\end{lemma}

\begin{proof} 
Recall that $\beta\in H^2$ is called a Seiberg-Witten basic class 
if $\exp(\beta)$ appears in the Seiberg-Witten invariant with nonzero
coefficient. Given this, the Seiberg-Witten basic classes of 
$X(d_1,d_2,d_3)$ are the classes $2\sum_{j=1}^3 b_jd_j[T_j]$ where
$b_j=-1,0$, or $1$. 

We first observe that for any orientation compatible symplectic structure 
$\omega$ on $X(d_1,d_2,d_3)$, the canonical class $c_1(K)$ must equal 
$2\sum_{j=1}^3 b_jd_j[T_j]$ where each of $b_1,b_2,b_3$ is nonzero. 
The reason is as follows. According to Taubes \cite{T2}, 
for any complex line bundle $E$, if $2 c_1(E)-c_1(K)$ is a Seiberg-Witten 
basic class, then the Poincar\'{e} dual of $c_1(E)$ is represented by the 
fundamental class of a symplectic submanifold; in particular, $c_1(E)\cdot 
[\omega]>0$ if $c_1(E)\neq 0$. Now observe that if say 
$c_1(K)=2(d_2[T_2]+d_3[T_3])$ (i.e. $b_1=0$), then since
both $2(\pm d_1[T_1]-d_2[T_2]-d_3[T_3])$ are Seiberg-Witten basic
classes, both $\pm d_1[T_1]$ have a positive cup product with
$[\omega]$, which is a contradiction. 

By replacing $[T_j]$ with $-[T_j]$ if necessary, we may assume without
loss of generality that $c_1(K)=2\sum_{j=1}^3 d_j[T_j]$. With this 
understood, note that for each $j=1,2,3$, $2d_j[T_j]-2\sum_{k=1}^3 d_k[T_k]$ 
is a Seiberg-Witten basic class, hence by Taubes' theorem in \cite{T2},
$d_j[T_j]$ is Poincar\'{e} dual to the fundamental class of a
symplectic submanifold, which implies that $[T_j]\cdot [\omega]>0$.
The lemma follows easily.

\end{proof}

\begin{lemma}
{\em(1)} 
Let $f:X(d_1,d_2,d_3)\rightarrow X(d_1^\prime,d_2^\prime,d_3^\prime)$
be any diffeomorphism. Then for $j=1,2,3$, one has $d_j=d_j^\prime$ and
$f^\ast ([T_j^\prime])=\pm [T_j]$. {\em(}Here $[T_j^\prime]$ denotes the
corresponding class of $X(d_1^\prime,d_2^\prime,d_3^\prime).${\em)}
In particular, $X(d_1,d_2,d_3)$ are distinct smooth $4$-manifolds
for distinct triples $(d_1,d_2,d_3)$.

{\em(2)} 
Let $\omega$ be any orientation compatible symplectic structure on 
$X(d_1,d_2,d_3)$ and $f$ be any self-diffeomorphism such that $f^\ast
[\omega]=[\omega]$. Then $f^\ast[T_j]=[T_j]$ for $j=1,2,3$.
\end{lemma}

\begin{proof}
Observe that $f$ must be orientation-preserving, because under an 
orientation-reversing diffeomorphism the signature changes by a sign of 
$-1$. Consequently,
$f^\ast$ sends the Seiberg-Witten basic classes of $X(d_1^\prime,
d_2^\prime,d_3^\prime)$ to those of $X(d_1,d_2,d_3)$. In particular,
there are $b_j\in\Z$, $j=1,2,3$, with each $b_j=-1,0$ or $1$ such that
$$
f^\ast(2d_1^\prime [T_1^\prime])=2\sum_{j=1}^3 b_jd_j[T_j].
$$
By Lemma 2.1 (1), there are homology classes $v_1,v_2,v_3$
such that $v_i\cdot [T_j]=1$ for $i=j$ and $v_i\cdot [T_j]=0$ otherwise.
Taking cup product of each side of the above equation with $v_1,v_2,v_3$, we
see that $d_1^\prime$ is a divisor of $d_j$ if $b_j\neq 0$. Since by
assumption $d_1^\prime>1$ and $d_1, d_2, d_3$ are pairwise relative prime, 
it follows that there exists exactly one $b_j$ which is 
nonzero. Applying the same argument to $f^{-1}$, we see that one actually has 
$d_1^\prime=d_j$ and $f^\ast([T_1^\prime])=\pm [T_j]$. Since each of the
triples $(d_1,d_2,d_3)$ and $(d_1^\prime,d_2^\prime,d_3^\prime)$ is assumed 
to be in the ascending order, we must have  
$d_j^\prime=d_j$ and $f^\ast([T_j^\prime])=\pm [T_j]$ for $j=1,2,3$, as 
claimed in (1).

If $\omega$ is an orientation compatible symplectic structure on 
$X(d_1,d_2,d_3)$ and $f$ is a self-diffeomorphism such that $f^\ast
[\omega]=[\omega]$, then $f^\ast[T_j]=[T_j]$, $j=1,2,3$, must be true
because $[T_j]\cdot [\omega]\neq 0$ by Lemma 2.2.

\end{proof}

In the remaining sections, we will abbreviate the notation and denote the
exotic $K3$ $X(d_1,d_2,d_3)$ by $X_\alpha$. 

\section{Recollection of various known results}

In this section we collect some theorems (known to date) and some 
observations that are scattered in the literature, which may be used 
to provide obstructions to the existence of 
certain smooth finite group actions on $4$-manifolds. (In fact, many of these 
obstructions also apply to locally linear topological actions.) For symplectic actions of 
a finite group on a minimal symplectic $4$-manifold with $c_1^2=0$, there are 
further results in terms of the fixed-point set structure of the action.
These will be briefly reviewed at the beginning of Section $6$, and details
may be found in \cite{CK}. 

\vspace{2mm}

{\it Borel spectral sequence.} \hspace{2mm} 
We review here some relevant results about locally linear topological 
actions of a finite group on a closed simply-connected $4$-manifold. The main 
technique for deriving these results is the Borel spectral sequence, cf.
e.g. \cite{Bre}. 

Let $G\equiv \Z_p$, where $p$ is prime, act locally linearly on a closed 
simply-connected $4$-manifold $M$ via orientation-preserving homeomorphisms, 
and let $F$ be the fixed-point set of the action. We first review a result 
due to Edmonds which describes a relationship between the fixed-point 
set $F$ and the existence of certain types of representations of $G$ on 
$H^2(M)$ induced by the action of $G$ on $M$. 

Recall that by a result of Kwasik and Schultz (cf. \cite{KS}), each 
integral representation of $\Z_p$ on $H^2(M)$ can be expressed as a sum of 
copies of the group ring $\Z[\Z_p]$ of $\Z$-rank $p$, the trivial 
representation $\Z$ of $\Z$-rank $1$, and the representation $\Z[\mu_p]$ 
of cyclotomic type of $\Z$-rank $p-1$, which is the kernel of the 
augmentation homomorphism $\Z[\Z_p]\rightarrow \Z$. Here $\mu_p\equiv 
\exp(\frac{2\pi i}{p})$, which will be used throughout. 

\begin{proposition} {\em(}cf. \cite{E}, Prop. 2.4{\em)} 
Assume that $F$ is nonempty. Let $b_1(F)$ be the first Betti number 
of $F$ in $\Z_p$-coefficients and let $c$ be the number of copies of 
$G$-representations of cyclotomic type in $H^2(M)$. Then $b_1(F)=c$. 
In particular, $c=0$ if the $G$-action is pseudofree, and $b_1(F)=0$ 
if the $G$-action is homologically trivial.
\end{proposition}

Another result of Edmonds gives some homological restrictions on the
$2$-dimensional components of the fixed-point set $F$.

\begin{proposition} {\em(}cf. \cite{E}, Cor. 2.6{\em)} 
If $F$ is not purely $2$-dimensional, then the $2$-dimensional components 
of $F$ represent independent elements of $H_2(M;\Z_p)$. If $F$ is purely 
$2$-dimensional, and has $k$ $2$-dimensional components, then the 
$2$-dimensional components span a subspace of $H_2(M;\Z_p)$ of dimension
at least $k-1$, with any $k-1$ components representing independent elements.
\end{proposition}

The next theorem, due to McCooey \cite{McC}, is concerned with locally
linear, homologically trivial topological actions by a compact Lie group
(e.g. a finite group) on a closed $4$-manifold.

\begin{theorem}
Let $G$ be a {\em(}possibly finite{\em)} 
compact Lie group, and suppose $M$ is a 
closed $4$-manifold with $H_1(M;\Z)=0$ and $b_2(M)\geq 2$, equipped with
an effective, locally linear, homologically trivial $G$-action. Denote
by $F$ the fixed-point set of $G$.
\begin{itemize}
\item [{1.}] If $b_2(M)=2$ and $F\neq\emptyset$, then $G$ is
isomorphic to a subgroup of $\s^1\times\s^1$.
\item [{2.}] If $b_2(M)\geq 3$, then $G$ is isomorphic to a subgroup of 
$\s^1\times\s^1$, and $F\neq\emptyset$.
\end{itemize}
\end{theorem}

\vspace{2mm}

{\it $G$-index theorems.} \hspace{2mm} Here we collect some formulas
which fall into the realm of $G$-index theorems of Atiyah and 
Singer (cf. \cite{AS}). In particular, these formulas allow us to relate 
the fixed-point set structure of the group action with the induced 
representation on the rational cohomology of the manifold.

Let $M$ be a closed, oriented smooth $4$-manifold, and let $G\equiv\Z_p$
be a cyclic group of prime order $p$ acting on $M$ effectively via
orientation-preserving diffeomorphisms. Then the fixed-point set $F$,
if nonempty, will be in general a disjoint union of finitely many isolated 
points and orientable surfaces. Fix a generator $g\in G$.
Then each isolated fixed point $m\in F$ is associated with a pair of 
integers $(a_m,b_m)$, where $0<a_m,b_m<p$, such that the action 
of $g$ on the tangent space at $m$ is given by the complex linear 
transformation $(z_1,z_2)\mapsto (\mu_p^{a_m}z_1,\mu_p^{b_m}z_2)$. (Note 
that $a_m,b_m$ are uniquely determined up to a change of order or a 
simultaneous change of sign modulo $p$.) Likewise, at each connected surface 
$Y\subset F$, there is an integer $c_Y$ with $0<c_Y<p$, 
which is uniquely determined up to a sign modulo $p$, such that the action of $g$ on 
the normal bundle of $Y$ in $M$ is given by $z\mapsto \mu_p^{c_Y} z$.

\begin{theorem} {\em(}Lefschetz Fixed Point Theorem{\em)}. 
$L(g,M)=\chi(F)$, where $\chi(F)$ is the Euler characteristic of the 
fixed-point set $F$ and $L(g,M)$ is the Lefschetz number of the map 
$g:M\rightarrow M$, which is defined by
$$
L(g,M)=\sum_{k=0}^4 (-1)^k tr(g)|_{H^k(M;\R)}.
$$
\end{theorem}

Note that the above theorem holds true for topological actions, cf. \cite{KS}.

\begin{theorem} {\em(}G-signature Theorem{\em)}. 
Set 
$$
\text{Sign}(g,M)=tr(g)|_{H^{2,+}(M;\R)}-tr(g)|_{H^{2,-}(M;\R)}.
$$
Then
$$
\text{Sign}(g,M)=\sum_{m\in F} -\cot(\frac{a_m\pi}{p})\cdot 
                      \cot(\frac{b_m\pi}{p})+ \sum_{Y\subset F}
                      \csc^2(\frac{c_Y\pi}{p})\cdot (Y\cdot Y),
$$
where $Y\cdot Y$ denotes the self-intersection number of $Y$.
\end{theorem}

Note that the G-signature Theorem is also valid for locally linear,
topological actions in dimension $4$, cf. e.g. \cite{Gor}. 

One can average the formula for $\text{Sign}(g,M)$ over $g\in G$ to
obtain the following version of the $G$-signature Theorem.

\begin{theorem}  {\em(}G-signature Theorem -- the weaker version{\em)}.  
$$
|G|\cdot\text{Sign}(M/G)=\text{Sign}(M)
+\sum_{m\in F}\text{def}_m +\sum_{Y\subset F}\text{def}_Y.
$$
where the terms $\text{def}_m$ and 
$\text{def}_Y$ {\em(}called  {\em signature defects)} are given 
by the following formulae:
$$
\text{def}_m=\sum_{k=1}^{p-1}\frac{(1+\mu_p^k)(1+\mu_p^{kq})}
{(1-\mu_p^k)(1-\mu_p^{kq})}
$$
if the local representation of $G$ at $m$ is given by $(z_1,z_2)\mapsto 
(\mu_p^{k}z_1,\mu_p^{kq}z_2)$, and 
$$
\text{def}_Y=\frac{p^2-1}{3}\cdot (Y\cdot Y).
$$
\end{theorem}

The above version of the G-signature Theorem is more often used because
the signature defect $\text{def}_m$ can be computed in terms of Dedekind
sum, cf. \cite{HZ}.

Now suppose that the $4$-manifold $M$ is spin, and that the $G$-action 
on $M$ lifts to the spin structures on $M$. Then the index of Dirac
operator $\D$ gives rise to a character of $G$. More precisely, for each 
$g\in G$, one can define the ``Spin-number'' of $g$ by
$$
\text{Spin}(g,M)=tr(g)|_{\text{Ker} \D}-tr(g)|_{\text{Coker} \D}.
$$
If we write $\text{Ker} \D=\oplus_{k=0}^{p-1} V_k^{+}$,
$\text{Coker} \D=\oplus_{k=0}^{p-1} V_k^{-}$, where $V_k^{+}$, $V_k^{-}$
are the eigenspaces of $g$ with eigenvalue $\mu_p^k$, then
$$
\text{Spin}(g,M)=\sum_{k=0}^{p-1} d_k\mu_p^k,
$$
where $d_k\equiv \dim_\C V_k^{+}-\dim_\C V_k^{-}$. Since both $\text{Ker} \D$
and $\text{Coker} \D$ are quaternion vector spaces, and the quaternions $i$
and $j$ are anti-commutative, it follows that $V_0^{\pm}$ are quaternion 
vector spaces, and when $p=|G|$ is odd, $j$ maps $V_k^{\pm}$ isomorphically 
to $V_{p-k}^{\pm}$ for $1\leq k\leq p-1$. This particularly implies that 
$d_0$ is even, and when $p$ is odd, $d_k=d_{p-k}$ for $1\leq k\leq p-1$.

\begin{theorem} {\em(}G-index Theorem for Dirac Operators, cf. \cite{AB}{\em)}.
Assume further that the action of $G$ on $M$ is spin and that there are only
isolated fixed points. Then the ``Spin-number'' 
$\text{Spin}(g,M)=\sum_{k=0}^{p-1} d_k\mu_p^k$ is given in terms of the 
fixed-point set structure by the following formula
$$
\text{Spin}(g,M)=-\sum_{m\in F}\epsilon(g,m)
                      \cdot \frac{1}{4}\csc(\frac{a_m\pi}{p})
                          \cdot\csc(\frac{b_m\pi}{p}),
$$
where $\epsilon(g,m)=\pm 1$ depends on the fixed point $m$ and the lifting
of the action of $g$ to the spin structure. 
\end{theorem} 

We give a formula below for the sign $\epsilon(g,m)$ with the assumption that
the action of $G$ preserves an almost complex structure on $M$ (e.g. the 
action of $G$ is via symplectic symmetries) and that the order of $G$ is odd.

\begin{lemma}
Assume further that $M$ is simply-connected, spin, 
and $|G|=p$ is an odd prime. Then the action of $G$ on $M$ is spin.
Moreover, if $M$ is almost complex and the action of $G$ preserves the 
almost complex structure on $M$, then the
``Spin-number'' can be computed by
\begin{eqnarray*}
\text{Spin}(g,M) & = & -\sum_{m\in F}\epsilon(g,m)
                      \cdot \frac{1}{4}\csc(\frac{a_m\pi}{p})
                          \cdot\csc(\frac{b_m\pi}{p})\\
                &   & +\sum_{Y\subset F}\epsilon(g,Y)\cdot\frac{(Y\cdot Y)}{4}
                 \csc(\frac{c_Y\pi}{p})\cdot \cot(\frac{c_Y\pi}{p}),
\end{eqnarray*}
where the signs $\epsilon(g,m)$ and $\epsilon(g,Y)$ are determined as
follows. First, in the above formula, $a_m,b_m$ and $c_Y$, which satisfy 
$0<a_m,b_m<p$ and $0<c_Y<p$, are uniquely determined because the 
corresponding complex representations 
$(z_1,z_2)\mapsto (\mu_p^{a_m}z_1,\mu_p^{b_m}z_2)$ and
$z\mapsto \mu_p^{c_Y} z$ are chosen to be compatible with the almost complex 
structure which $G$ preserves. With this convention, $\epsilon(g,m)$ and 
$\epsilon(g,Y)$ are given by
$$
\epsilon(g,m)= (-1)^{k(g,m)},\hspace{3mm} 
\epsilon(g,Y)= (-1)^{k(g,Y)}
$$
where $k(g,m)$ and $k(g,Y)$ are defined by equations 
$$
k(g,m)\cdot p=2r_m+a_m+b_m, \hspace{3mm} 
k(g,Y)\cdot p=2r_Y+c_Y
$$
for some $r_m$, $r_Y$ satisfying $0\leq r_m<p$ and $0<r_Y<p$. 
\end{lemma}

\begin{proof}
We first show that the action of $G$ is spin. Let $E_G\rightarrow B_G$
be the universal principal $G$-bundle. Then observe that a bundle $E$ 
over $M$ as a $G$-bundle corresponds to a bundle $E^\prime$ over
$E_G\times_G M$ whose restriction to the fiber $M$ of the fiber bundle
$E_G\times_G M\rightarrow B_G$ is $E$. With this understood, a $G$-spin
structure on $M$ corresponds to a principal $Spin(4)$-bundle over
$E_G\times_G M$ whose restriction to the fiber $M$ is a spin structure
on $M$. The obstruction to the existence of such a bundle is a class
in $H^2(E_G\times_G M;\Z_2)$ which maps to the second Stiefel-Whitney
class $w_2(M)\in H^2(M;\Z_2)$ under the homomorphism
$i^\ast: H^2(E_G\times_G M;\Z_2)\rightarrow H^2(M;\Z_2)$ induced by 
the inclusion $i:M\rightarrow E_G\times_G M$ as a fiber. The obstruction 
vanishes because (1) $w_2(M)=0$ since $M$ is spin, (2) the homomorphism 
$i^\ast:H^\ast(E_G\times_G M;\Z_2)\rightarrow H^\ast(M;\Z_2)$ is a
monomorphism, which can be seen easily by the transfer argument (cf. 
\cite{Bre}) given that the order $|G|=p$ is odd. This proves 
that the action of $G$ is spin. Note that there is a unique $G$-spin 
structure because the $G$-spin structures are classified by 
$H^1(E_G\times_G M;\Z_2)=0$. Here $H^1(E_G\times_G M;\Z_2)=0$ because 
$H^1(M;\Z_2)=0$ (since $M$ is simply-connected) and 
$i^\ast:H^1(E_G\times_G M;\Z_2)\rightarrow H^1(M;\Z_2)$ is injective.

Now suppose $M$ is almost complex and $G$ preserves the almost complex 
structure on $M$. Then the 
$G$-spin structure as the unique $G$-$Spin^\C$ structure with trivial
determinant line bundle is given by a (unique) $G$-complex line bundle
$L$ over $M$ such that $L^2\otimes K^{-1}=M\times\C$ as a $G$-bundle where
the action of $G$ on $\C$ is trivial. Here $K$ is the canonical bundle
of the almost complex structure. Moreover, the Dirac operator $\D$ is
simply given by the $\bar{\partial}$-complex twisted with the complex line
bundle $L$. The ``Spin-number'' $\text{Spin}(g,M)$ may be computed using
the $G$-index Theorem for the $\bar{\partial}$-complex (i.e. the 
holomorphic Lefschetz fixed point theorem), cf. \cite{AS}. 

More concretely, let the action of $g$ at a fixed point $m\in F$ and a
fixed component $Y\subset F$ be denoted by $z\mapsto \mu_p^{r_m}z $
and $z\mapsto \mu_p^{r_Y}z $ respectively. Then $L^2\otimes K^{-1}=M\times\C$
as a trivial $G$-bundle implies that 
$$
2r_m+a_m+b_m=0\pmod{p}, \hspace{3mm} 2r_Y+c_Y=0\pmod{p}.
$$
We shall impose further conditions that $0\leq r_m<p$ and $0<r_Y<p$, and 
define integers $k(g,m)$, $k(g,Y)$ as in the lemma by
$$
k(g,m)\cdot p=2r_m+a_m+b_m, \hspace{3mm} 
k(g,Y)\cdot p=2r_Y+c_Y.
$$
With these understood, the contribution to $\text{Spin}(g,M)$ from $m\in F$ is
$$
I_m=\frac{\mu_p^{r_m}}
{(1-\mu_p^{-a_m})(1-\mu_p^{-b_m})}=(-1)^{k(g,m)+1}\cdot
                          \frac{1}{4}\csc(\frac{a_m\pi}{p})
                          \cdot\csc(\frac{b_m\pi}{p}),
$$
and the contribution from $Y\subset F$ is 
$$
I_Y=\frac{\mu_p^{r_Y}(1+l)(1+t/2)}{1-\mu_p^{-c_Y}(1-n)}[Y]
               =(-1)^{k(g,Y)}\cdot\frac{(Y\cdot Y)}{4}
                 \csc(\frac{c_Y\pi}{p})\cdot \cot(\frac{c_Y\pi}{p}),
$$
where $l,t,n$ are the first Chern classes of $L$, $TY$ and the normal
bundle of $Y$ in $M$. The formula for $\text{Spin}(g,M)$ follows 
immediately.

\end{proof}

\vspace{2mm}

{\it Seiberg-Witten equations.} \hspace{2mm}
There are obstructions to the existence of smooth finite group actions on
$4$-manifolds that come from Seiberg-Witten theory, based on the ideas in
Furuta \cite{Fu}. See \cite{Bry, Fa, FF, Na}. 

\begin{theorem}{\em (}cf. \cite{Fa, Na}{\em)}\hspace{2mm}
Let $M$ be a closed, oriented smooth $4$-manifold with $b_1=0$ and
$b_2^{+}\geq 2$, which admits a smooth $G\equiv\Z_p$ action of prime 
order such that $H^{2}(M;\R)$ contains a $b_2^{+}$-dimensional subspace
consisting of invariant elements of positive square. Let $c$ be a 
$G$-$Spin^\C$ structure on $M$ such that the $G$-index of the Dirac 
operator $\text{ind}_G\D=\sum_{k=0}^{p-1}d_k\C_k$ satisfies 
$2d_k\leq b_2^{+}-1$ for all $0\leq k<p$.
Then the corresponding Seiberg-Witten invariant obeys
$$
\text{SW}_M(c)=0 \pmod{p}.
$$
Here $\C_k$ denotes the complex $1$-dimensional weight $k$ representation
of $G\equiv\Z_p$. 
\end{theorem}

\begin{theorem}{\em (}cf. \cite{Fu, FF}{\em)}\hspace{2mm}
Suppose a smooth action of a finite group $G$ on a closed, spin $4$-manifold 
$M$ is spin. Let $\D$ be the Dirac operator on the spin $4$-orbifold $M/G$. 
Then either
$\text{ind}\;\D=0$ or $-b_2^{-}(M/G)<\text{ind}\;\D<b_2^{+}(M/G)$.
\end{theorem}

We remark that when the action of $G$ preserves an almost complex 
structure on $M$, the index of the Dirac operator, $\text{ind}\;\D$,
for the $4$-orbifold $M/G$ can be calculated using Lemma 3.8, or using
the formula for the dimension of the corresponding Seiberg-Witten moduli 
space in \cite{C0}. (See also \cite{C1}.)

\vspace{2mm}

{\it The Kirby-Siebenmann and the Rochlin invariants.}\hspace{2mm} Suppose 
a locally linear, topological action of a finite group $G$ on a $4$-manifold $M$ 
is spin and pseudofree. Then the quotient space $M/G$ is a spin $4$-orbifold
with only isolated singular points. Let $N$ be the spin $4$-manifold with boundary
obtained from $M/G$ by removing a regular neighborhood of the singular set, and
denote by $\partial\eta$ the spin structure on $\partial N$ induced from that of $N$.
Then the Kirby-Siebenmann invariant of $N$, denoted by $\text{ks}(N)$, and the 
Rochlin invariant of $(\partial N,\partial\eta)$, denoted by $\text{roc}(\partial N,\partial\eta)$, 
are constrained by the following 

\begin{theorem}{\em (}cf. \S 10.2B in \cite{FQ}{\em)}\hspace{2mm}
$$
8\cdot\text{ks}(N)\equiv \text{Sign}(N) +\text{roc}(\partial N,\partial\eta) \pmod{16}.
$$
\end{theorem}

Note that a necessary condition for the $G$-action to be smoothable is $\text{ks}(N)=0$.

\section{Smooth cyclic actions}

In this section, we give proofs of Theorem 1.1 and Theorem 1.4. 

\vspace{3mm}

The following lemma, together with Lemma 2.3 (1), settles Theorem 1.1 (1)
because Lemma 2.3 (1) implies that the classes $[T_j]$, $j=1,2,3$, are fixed 
under the action of  $\Z_p$ whenever $p$ is odd. 

\begin{lemma}
For any smooth action of a finite group $G$ on $X_\alpha$ which fixes the 
classes $[T_j]$, $j=1,2,3$, there is a $3$-dimensional subspace of 
$H^{2}(X_\alpha;\R)$ which is fixed under $G$ and over which the cup-product
is positive definite. 
\end{lemma}

\begin{proof}
By Lemma 2.1 (1), and since $H^2(X_\alpha;\Z)$ and $H^2(X;\Z)$ are naturally
identified, there are classes $v_i$, $i=1,2,3$, in $H^2(X_\alpha;\Z)$ such 
that $v_i\cdot [T_j]=1$ if $i=j$ and $v_i\cdot [T_j]=0$ otherwise. 

For any given $g\in G$, we set 
$v_j^\prime\equiv\sum_{k=0}^{|g|-1}(g^k)^\ast v_j$. 
Then $g^\ast v_j^\prime=v_j^\prime$. Now for sufficiently small $\epsilon>0$, 
we obtain three linearly independent classes
$$
u_j\equiv [T_j]+\epsilon v_j^\prime, \; j=1,2,3,
$$
which are all fixed by $g$, i.e., $g^\ast u_j=u_j$ for $j=1,2,3$. 
On the other hand, 
set $c\equiv \max\{|(v_1^\prime)^2|+1,|(v_2^\prime)^2|+1,
|(v_3^\prime)^2|+1\}>0$. Then for any $i,j$,
\begin{eqnarray*}
u_i\cdot u_j & = & [T_i]\cdot [T_j]+\epsilon ([T_i]\cdot v_j^\prime+
[T_j]\cdot v_i^\prime)+\epsilon^2 v_i^\prime\cdot v_j^\prime\\
             & = & \epsilon |g|([T_i]\cdot v_j+[T_j]\cdot v_i)
                   +\epsilon^2 v_i^\prime\cdot v_j^\prime\\
             &= & \left \{\begin{array}{cc}
2\epsilon|g|+\epsilon^2 v_i^\prime\cdot v_j^\prime & \mbox{ if } i=j\\
\epsilon^2 v_i^\prime\cdot v_j^\prime  & \mbox{ if } i\neq j\\
\end{array} \right .
\end{eqnarray*}
and for any $a_1, a_2, a_3\in \R$ and any $0<\epsilon <10^{-1}c^{-1}$,
\begin{eqnarray*}
(\sum_{i=1}^3 a_i u_i)^2 & \geq & 2\epsilon|g|(\sum_{i=1}^3 a_i^2)
-3c\cdot \epsilon^2(\sum_{i=1}^3 a_i^2)\\
                         & \geq & \epsilon (\sum_{i=1}^3 a_i^2).\\
\end{eqnarray*}
Consequently, for sufficiently small $\epsilon>0$, $u_1,u_2,u_3$ 
span a $3$-dimensional subspace of $H^{2}(X_\alpha;\R)$ which is fixed 
under the action of $G$ and over which the cup-product is positive
definite.

\end{proof}

Recall that the Kummer surface $X$ (as well as the exotic $X_\alpha$) 
has intersection form $3H\oplus 2(-E_8)$, where $H=\left (
\begin{array}{ll}
0 & 1\\
1 & 0
\end{array} \right )$ and

\vspace{3mm}

$$
-E_8=\left (\begin{array}{llllllll}
-2 & 1 & & & & & & \\
1 & -2 & 1 & & & & &\\
 & 1 & -2 & 1 & & & & \\
 & & 1 & -2 & 1 & & & \\
 & & & 1 & -2 & 1 & & 1\\
 & & & & 1 & -2 & 1 & \\
 & & & & & 1 & -2 & \\
 & & & & 1 & & & -2
\end{array} \right ).
$$

\vspace{3mm}

The $-E_8$ form is the intersection matrix of a standard basis 
$\{f_j|1\leq j\leq 8\}$, which may be conveniently described by the
graph in Figure 1, where two nodes are connected by an edge if and only
if the corresponding basis vectors have intersection product $1$.

\begin{figure}[ht]

\setlength{\unitlength}{1pt}

\begin{picture}(420,100)(-210,-195)

\linethickness{0.5pt}

\put(-120,-140){$f_1$}

\put(-80,-140){$f_2$}

\put(-40,-140){$f_3$}

\put(0,-140){$f_4$}

\put(40,-140){$f_5$}

\put(80,-140){$f_6$}

\put(120,-140){$f_7$}

\put(50,-190){$f_8$}

\put(-120,-150){\circle*{3}}

\put(-120,-150){\line(1,0){240}}

\put(-80,-150){\circle*{3}}

\put(-40,-150){\circle*{3}}

\put(0,-150){\circle*{3}}

\put(40,-150){\circle*{3}}

\put(80,-150){\circle*{3}}

\put(120,-150){\circle*{3}}

\put(40,-150){\line(0,-1){40}}

\put(40,-190){\circle*{3}}

\end{picture}

\caption{}

\end{figure}

We shall next exhibit some geometric representative of a standard basis 
for each of the two $(-E_8)$-summands in the intersection form of the Kummer
surface $X$, which plays a crucial role in analyzing the induced action on 
$H_2(X_\alpha;\Z)$ of a smooth finite group action on the exotic $K3$ surface
$X_\alpha$.

\begin{lemma}
There exist two disjoint geometric representatives of a standard basis
of the $-E_8$ form in the Kummer surface $X$, which both lie in the 
complement of the tori $T_1$, $T_2$, $T_3$.
\end{lemma}

\begin{proof}
Recall that $X$ is the $4$-manifold obtained by replacing each
of the $16$ singularities $(\pm 1,\pm 1,\pm 1,\pm 1)$ in
$T^4/\langle \hat{\rho}\rangle$ by a $(-2)$-sphere.
We denote the $(-2)$-spheres by $\Sigma(\pm 1,\pm 1,\pm 1,\pm 1)$
accordingly and call them the exceptional $(-2)$-spheres in $X$.

On the other hand, recall from Section $2$ that for $j=1,2,3$, there
is a minimal complex surface $X(j)$ and an elliptic fibration 
$X(j)\rightarrow \s^2$, where $X(j)$ is obtained by resolving the
singularities of $T^4/\langle\hat{\rho}\rangle$ and the
elliptic fibration comes from the fibration
$\pi_j:T^4/\langle\hat{\rho}\rangle\rightarrow T^2/\langle\hat{\rho}\rangle$ 
induced by the projection
$$
(z_0,z_1,z_2,z_3)\mapsto (z_0,z_j).
$$
Note that $\pi_j:T^4/\langle\hat{\rho}\rangle\rightarrow T^2/\langle\hat{\rho}\rangle$ 
has $4$ singular fibers, which are over $(\pm 1,\pm 1)\in T^2/\langle\hat{\rho}\rangle$. 
We denote the proper transform of $\pi_j^{-1}(\pm 1,\pm 1)$ in $X(j)$ by 
$\Sigma_j (\pm 1,\pm 1)$, which is also a $(-2)$-sphere.

Recall also that for each $j$ we have fixed an identification between $X$ and 
the complex surface $X(j)$. Note that under such an identification each
exceptional $(-2)$-sphere in $X$ inherits an orientation from the 
corresponding complex curve in $X(j)$. For the purpose here we shall  
arrange the identifications between $X$ and the complex surfaces $X(j)$
such that each of the exceptional $(-2)$-spheres in $X$ inherits a
consistant orientation, and as a result, each of them is oriented and
defines a homology class in $H_2(X;\Z)$. With such identifications between 
$X$ and the complex surfaces $X(j)$ fixed, we shall regard the 
$(-2)$-spheres $\Sigma_j (\pm 1,\pm 1)$ in $X(j)$ as smooth surfaces 
in $X$, and call these $(-2)$-spheres the proper transform $(-2)$-spheres 
in $X$. We orient each $\Sigma_j (\pm 1,\pm 1)$ in $X$ by the
canonical orientation of the corresponding complex curve in $X(j)$.

With the choice of orientations on each $(-2)$-sphere (exceptional or
proper transform) understood, we observe
that (1) any two distinct exceptional $(-2)$-spheres have intersection 
product $0$ because they are disjoint, and (2) a proper transform 
$(-2)$-sphere and an exceptional $(-2)$-sphere have intersection product
either $0$ or $1$, depending on whether they are disjoint or not. 
The intersection product of two distinct proper transform $(-2)$-spheres
are described below.

\vspace{2mm}

{\bf Claim:} {\em Let $\kappa,\tau,\kappa^\prime,\tau^\prime$ 
take values in $\{1,-1\}$. Then the
following hold true: {\em(1)} If $(\kappa,\tau)\neq 
(\kappa^\prime,\tau^\prime)$, then $\Sigma_j(\kappa,\tau)$ and
$\Sigma_j(\kappa^\prime,\tau^\prime)$ are disjoint so that their intersection
product is $0$, {\em (2)} If $j\neq j^\prime$, then the intersection product 
of $\Sigma_j(\kappa,\tau)$ and $\Sigma_{j^\prime}(\kappa^\prime,\tau^\prime)$ 
is $0$ when $\kappa\neq\kappa^\prime$ {\em(}in fact the two $(-2)$-spheres are 
disjoint{\em)}, and is $-1$ when $\kappa=\kappa^\prime$.}

\vspace{2mm}

Accepting the above claim momentarily, one can easily verify that the 
following are two disjoint geometric representatives of a standard basis 
of the $-E_8$ form and that both lie in the complement of the three tori 
$T_1, T_2$ and $T_3$:
\begin{itemize}
\item [{(1)}] $f_1=-\Sigma_3(1,-1)-\Sigma(1,-1,-1,-1)-\Sigma(1,1,-1,-1)$, 
$f_2=\Sigma(1,1,-1,-1)$, $f_3=\Sigma_2(1,-1)+\Sigma(1,-1,-1,-1)$,
$f_4=\Sigma(1,1,-1,1)$, $f_5=\Sigma_3(1,1)+\Sigma(1,-1,-1,1)$, 
$f_6=\Sigma(1,1,1,1)$, $f_7=-\Sigma_2(1,1)-\Sigma(1,1,1,-1)
-\Sigma(1,1,1,1)$, $f_8=\Sigma_1(1,-1)+\Sigma(1,-1,1,1)$
\item [{(2)}] $f_1=-\Sigma_3(-1,-1)-\Sigma(-1,-1,-1,-1)-\Sigma(-1,1,-1,-1)$, 
$f_2=\Sigma(-1,1,-1,-1)$, $f_3=\Sigma_2(-1,-1)+\Sigma(-1,-1,-1,-1)$,
$f_4=\Sigma(-1,1,-1,1)$, $f_5=\Sigma_3(-1,1)+\Sigma(-1,-1,-1,1)$, 
$f_6=\Sigma(-1,1,1,1)$, $f_7=-\Sigma_2(-1,1)-\Sigma(-1,1,1,-1)
-\Sigma(-1,1,1,1)$, $f_8=\Sigma_1(-1,-1)+\Sigma(-1,-1,1,1)$
\end{itemize}

It remains to verify the claim. Note that part (1) of the claim follows from
the fact that the two proper transform $(-2)$-spheres lie in two distinct
fibers of the $C^\infty$-elliptic fibration $\pi_j:X\rightarrow \s^2$.
To see part (2), we suppose $j\neq j^\prime$. Then $\Sigma_j(\kappa,\tau)$ 
and $\Sigma_{j^\prime}(\kappa^\prime,\tau^\prime)$ are disjoint if $\kappa
\neq \kappa^\prime$ (because $\kappa$, $\kappa^\prime$ are the 
$z_0$-coordinates), and part (2) of the claim holds true in this case. 
Therefore we shall assume $\kappa=\kappa^\prime$. Without loss of generality, 
we may assume that $\kappa=\kappa^\prime=1$, and for simplicity we shall 
only check the case where $\tau=\tau^\prime=1$ and $j=2$, $j^\prime=3$. 
With this understood, note that the fiber class of $\pi_2:X\rightarrow\s^2$,
which is the class of the torus $T_2$, equals 
$$
2\cdot\Sigma_2(1,1)+\Sigma(1,1,1,1)+\Sigma(1,1,1,-1)
+\Sigma(1,-1,1,1)+\Sigma(1,-1,1,-1)
$$
and the fiber class of $\pi_3:X\rightarrow\s^2$,
which is the class of the torus $T_3$, equals 
$$
2\cdot\Sigma_3(1,1)+\Sigma(1,1,1,1)+\Sigma(1,1,-1,1)
+\Sigma(1,-1,1,1)+\Sigma(1,-1,-1,1).
$$
(Note that each $C^\infty$-elliptic fibration $\pi_j:X\rightarrow \s^2$ 
has $4$ singular fibers, all of type $I_0^\ast$, cf. \cite{BPV}, page 201.)
The assertion that the intersection product of $\Sigma_2(1,1)$ and 
$\Sigma_3(1,1)$ equals $-1$ follows immediately from the fact that
$[T_2]\cdot [T_3]=0$. This finishes the verification of the claim above,
and the proof of Lemma 4.2 is completed. 

\end{proof}

\begin{lemma}
Let $G$ be a finite group acting on $H_2(X_\alpha;\Z)$ preserving the 
intersection form and fixing each $[T_j]$, $j=1,2,3$. Then there is an
induced homomorphism $\Theta:G\rightarrow \mbox{Aut }(E_8\oplus E_8)$
such that the action of $G$ on $H_2(X_\alpha;\Z)$ is trivial if and only if
the induced homomorphism $\Theta$ has trivial image.
\end{lemma}

\begin{proof}
Let $\xi_k, \eta_k$, $1\leq k\leq 8$, be the classes in $H_2(X_\alpha;\Z)$ 
corresponding to the two standard bases of the $-E_8$ form defined in the 
previous lemma. Then the intersection form 
of $X_\alpha$ is isomorphic to $3H$ when restricted to the orthogonal 
complement of $\mbox{Span }(\xi_k,\eta_k|1\leq k\leq 8)$. By Lemma 4.2 and
Lemma 2.1 (1), there are classes $w_i\in H_2(X_\alpha;\Z)$, $i=1,2,3$, 
such that
$$
w_i\cdot [T_j]=\left \{\begin{array}{ll}
0 & \mbox{if } i\neq j\\
1 & \mbox{if } i=j,
\end{array}\right.\;\; w_i\cdot \xi_k=w_i\cdot \eta_k=0 \mbox{ and } 
w_i\cdot w_j=0 \;\; \forall i, j,k.
$$

Let $\Omega$ be the orthogonal complement of 
$\mbox{Span}([T_1],[T_2],[T_3])$. We shall prove that 
$\Omega=\mbox{Span }([T_j],\xi_k,\eta_k|j=1,2,3, 1\leq k\leq 8)$. 
To see this, observe that $w_1,w_2,w_3$, $[T_1]$, $[T_2]$, $[T_3]$,
and $\xi_k, \eta_k$, $1\leq k\leq 8$, form a basis of $H_2(X_\alpha;\Z)$.
For any class $x\in H_2(X_\alpha;\Z)$, expand $x$ in the above basis. 
Then by Lemma 4.2, there are no terms of $w_1,w_2,w_3$ in the expansion
of $x$ if and only if its intersection product with each of  $[T_1]$, 
$[T_2]$, $[T_3]$ is zero. This proves our claim about $\Omega$.

Since $G$ fixes each $[T_j]$, $j=1,2,3$, the orthogonal complement of 
$\mbox{Span}([T_1],[T_2],[T_3])$, which is 
$\mbox{Span }([T_j],\xi_k,\eta_k|j=1,2,3, 1\leq k\leq 8)$, is invariant
under the action of $G$. The induced action of $G$ on
$$
\mbox{Span }([T_j],\xi_k,\eta_k|j=1,2,3, 1\leq k\leq 8)/
\mbox{Span }([T_1],[T_2],[T_3])
$$ 
gives rise to a homomorphism 
$\Theta: G\rightarrow \mbox{Aut }(E_8\oplus E_8)$ by Lemma 4.2.

It remains to show that if $\Theta$ has trivial image, then the
action of $G$ on $H_2(X_\alpha;\Z)$ is also trivial. To see this,
let $g\in G$ be any given element. Then for each $\alpha\in
\mbox{Span }(\xi_k,\eta_k|1\leq k\leq 8)$, there exists a $u_\alpha\in
\mbox{Span }([T_1],[T_2],[T_3])$ such that
$g\cdot\alpha=\alpha+u_\alpha$. This gives, for each $n\in\Z^{+}$,
$$
g^n\cdot\alpha=\alpha +n u_\alpha \mbox{    because } g\cdot u_\alpha
=u_\alpha.
$$
It follows easily that $u_\alpha=0$ since $g$ is of finite order,
and consequently, $g\cdot \alpha=\alpha$.

For each $w_j$, $j=1,2,3$, there are $\hat{w}_j\in
\mbox{Span }(w_1,w_2,w_3)$, $u_j\in\mbox{Span }([T_1],[T_2],[T_3])$ 
and $\alpha_j\in \mbox{Span }(\xi_k,\eta_k|1\leq k\leq 8)$, such that
$$
g\cdot w_j=\hat{w}_j+u_j+\alpha_j.
$$
Taking intersection product with $[T_k]$, $k=1,2,3$, we see that
$\hat{w}_j=w_j$, and taking intersection product with $\alpha_j
=g\cdot \alpha_j$, we see that $\alpha_j^2=0$, hence $\alpha_j=0$ 
(because $-E_8$ is negative definite). Consequently, we have 
$g\cdot w_j=w_j+u_j$. Since $g$ fixes $u_j$ and is of finite order,
we see as in the earlier argument that $u_j=0$, and therefore 
$g\cdot w_j=w_j$. Thus $G$ acts trivially on $H_2(X_\alpha;\Z)$ if 
$\Theta$ has trivial image.

\end{proof}

The following corollary gives Theorem 1.1 (2) for the case of $p>7$.

\begin{corollary}
Let $G$ be a $p$-group with $p>7$. Then any smooth $G$-action on 
$X_\alpha$ must act trivially on homology; in particular, $G$ must be 
abelian of rank at most $2$. 
\end{corollary}

\begin{proof}
Since $|G|$ is odd, the classes $[T_j]$, $j=1,2,3$, are fixed
under the $G$-action by Lemma 2.3 (1). On the other hand, $\Theta:G
\rightarrow \text{Aut }(E_8\oplus E_8)$ must have trivial image because
the order of $\text{Aut }(E_8\oplus E_8)$ is 
$2^{29}\cdot 3^{10}\cdot 5^4\cdot 7^2$ (cf. \cite{Se}), which is not 
divisible by any prime $p>7$.  By Lemma 4.3, the induced $G$-action on
$H_2(X_\alpha;\Z)$ must be trivial. It follows that the $G$-action is 
homologically trivial because $H_1(X_\alpha;\Z)=H_3(X_\alpha;\Z)=0$
($X_\alpha$ is simply-connected). The last assertion follows from 
McCooey's theorem (cf. Theorem 3.3).

\end{proof}

The following lemma can be found in \cite{E1}, however, for completeness
we sketch its proof here.

\begin{lemma}
The following are the only possibilities for integral representations of
$\Z_p$ for $p=3,5,7$ induced by $\Z_p\subset \text{Aut }(E_8)$:
\begin{itemize}
\item [{$\Z_3$}]: $\Z[\Z_3]\oplus \Z^5$,  $\Z[\Z_3]^2\oplus \Z^2$, 
                  $\Z[\Z_3]\oplus\Z[\mu_3]^2\oplus\Z$, and $\Z[\mu_3]^4$
\item [{$\Z_5$}]: $\Z[\Z_5]\oplus \Z^3$ and $\Z[\mu_5]^2$
\item [{$\Z_7$}]: $\Z[\Z_7]\oplus \Z$.
\end{itemize}
\end{lemma}

\begin{proof}
Since $p<23$, by a result of Reiner (cf. \cite{CR}) such a representation 
is of the form $\Z[\Z_p]^r\oplus\Z[\mu_p]^s\oplus\Z^t$, where $pr+(p-1)s+t=8$.
By Hambleton and Riehm \cite{HR}, $s$ must be even. Moreover, observe that
$\Z[\mu_p]^s$ and $\Z^t$ are always orthogonal to each other, so that if $r=0$
one of $s$ or $t$ must be $0$ as well because the form $E_8$ is not
splittable. The lemma follows.

\end{proof}

We remark that the integral representations of $\Z_p$ in the above lemma 
are all realized by a subgroup of $\text{Aut }(E_8)$ of order $p$, cf.
\cite{E1}.  

\vspace{2mm}

The following proposition settles the case of $p=7$ in Theorem 1.1 (2).

\begin{proposition}
Suppose $G\equiv\Z_p$, where $p=3,5,7$, acts smoothly on $X_\alpha$
such that the integral $G$-representation given by $\Theta:G\rightarrow
\text{Aut}(E_8\oplus E_8)$ contains no summands of cyclotomic type.
Then the intersection form on $H_2(X_\alpha;\Z)$ may be decomposed 
as $3H\oplus 2(-E_8)$ such that each summand 
$H$ or $-E_8$ is invariant under the $G$-action. Moreover, the action 
of $G$ on each $H$-summand is trivial.
\end{proposition}

\begin{proof}
Recall that $\text{Aut }(E_8\oplus E_8)$
is a semi-direct product of $\text{Aut }(E_8)\times\text{Aut }(E_8)$
by $\Z_2$ (cf. \cite{Se}). Since the order $|G|=p$ is odd, $G$ maps trivially 
to $\Z_2$ under $\Theta:G\rightarrow \text{Aut }(E_8\oplus E_8)$
and hence it can not exchange the two $E_8$-summands. It follows that each of
$$
\mbox{Span }([T_j],\xi_k|j=1,2,3, 1\leq k\leq 8),\;
\mbox{Span }([T_j],\eta_k|j=1,2,3, 1\leq k\leq 8)
$$
is invariant under the action of $G$, and there are two induced integral 
representations of $G$ on $E_8$ given by the action on 
$$
\mbox{Span }([T_j],\xi_k|j=1,2,3, 1\leq k\leq 8)
/\mbox{Span }([T_1],[T_2],[T_3])
$$
and 
$$
\mbox{Span }([T_j],\eta_k|j=1,2,3, 1\leq k\leq 8)
/\mbox{Span }([T_1],[T_2],[T_3])
$$
respectively.

We claim that there are classes 
$\xi_k^\prime,\eta_k^\prime\in H_2(X_\alpha;\Z)$ 
such that 
\begin{itemize}
\item [{(i)}] $\xi_k^\prime=\xi_k, \eta_k^\prime=\eta_k\pmod{\text{Span }
([T_1],[T_2],[T_3])}$,
\item [{(ii)}] $\text{Span }(\xi_k^\prime| 1\leq k\leq 8)$, 
$\text{Span }(\eta_k^\prime| 1\leq k\leq 8)$ are invariant under $G$.
\end{itemize}
Note that $\text{Span }(\xi_k^\prime| 1\leq k\leq 8)$ and 
$\text{Span }(\eta_k^\prime| 1\leq k\leq 8)$ split off two $G$-invariant
copies of $-E_8$ from $H_2(X_\alpha;\Z)$. The orthogonal complement, which 
is isomorphic to $3H$ and is also $G$-invariant, contains 
$\text{Span }([T_1],[T_2],[T_3])$. A similar argument as in the proof of 
Lemma 4.3 shows that the action of $G$ is trivial on each copy of $H$. 

It remains to verify the above claim. For simplicity, we shall only consider 
the case of $\xi_k$'s, the other case is completely parallel. Let $g\in G$ be
a fixed generator.

The key point of the proof is that a summand of type $\Z$ or $\Z[\Z_p]$ in 
$$
\mbox{Span }([T_j],\xi_k|j=1,2,3, 1\leq k\leq 8)
/\mbox{Span }([T_1],[T_2],[T_3])
$$ 
can be lifted to a $\Z[\Z_p]$-submodule of the same type in 
$\mbox{Span }([T_j],\xi_k|j=1,2,3, 1\leq k\leq 8)$. 
By Lemma 4.5 these are the only types of summands if there are no summands 
of cyclotomic type (which is always the case when $p=7$). 

More concretely, let $x$ be a generator of a $\Z$-summand in
$$
\mbox{Span }([T_j],\xi_k|j=1,2,3, 1\leq k\leq 8)
/\mbox{Span }([T_1],[T_2],[T_3])
$$
and let $x^\prime$ be any lift of $x$ in 
$\mbox{Span }([T_j],\xi_k|j=1,2,3, 1\leq k\leq 8)$. Then 
$g\cdot x^\prime=x^\prime+u$ for some $u\in\mbox{Span }([T_1],[T_2],[T_3])$.
As we argued in the proof of Lemma 4.3, this implies that $u=0$ and 
$g\cdot x^\prime=x^\prime$. Hence $x^\prime$ generates a 
$\Z[\Z_p]$-submodule of the same type which 
is a lift of the original $\Z$-summand.

Let $y$ be a generator of a $\Z[\Z_p]$-summand (as a $\Z[\Z_p]$-submodule).
Pick any lift $y^\prime$ of $y$ in 
$\mbox{Span }([T_j],\xi_k|j=1,2,3, 1\leq k\leq 8)$, then $y^\prime$ generates 
a free $\Z[\Z_p]$-submodule in 
$\mbox{Span }([T_j],\xi_k|j=1,2,3, 1\leq k\leq 8)$
which is a lift of the original $\Z[\Z_p]$-summand of the same type. 

Now suppose the integral $G$-representation  
$$
\mbox{Span }([T_j],\xi_k|j=1,2,3, 1\leq k\leq 8)
/\mbox{Span }([T_1],[T_2],[T_3])
$$
is decomposed as $\Z^t\oplus\Z[\Z_p]^r$, and let
$\{x_i,y_j|1\leq i\leq t, 1\leq j\leq r\}$
be a set of generators of the summands as $\Z[\Z_p]$-submodules. Then the set
$$
\{x_i,y_j,g\cdot y_j,\cdots,g^{p-1}\cdot y_j\}
$$
forms a $\Z$-basis of
$$
\mbox{Span }([T_j],\xi_k|j=1,2,3, 1\leq k\leq 8)
/\mbox{Span }([T_1],[T_2],[T_3]).
$$
Note that the intersection form on 
$$
\text{Span }(x_i^\prime,y_j^\prime,g\cdot y_j^\prime,
\cdots,g^{p-1}\cdot y_j^\prime),
$$
i.e., the span of the lifts, is isomorphic to that on 
$$
\mbox{Span }([T_j],\xi_k|j=1,2,3, 1\leq k\leq 8)
/\mbox{Span }([T_1],[T_2],[T_3]).
$$
The existence of $\xi_k^\prime$'s follows immediately.

This completes the proof of the proposition.

\end{proof}

\begin{remark}
In general, a summand of cyclotomic type may not be lifted to a summand
of the same type under a quotient homomorphism. For a simple example,
let us consider the integral $\Z_2$-representation on 
$\Z\langle x\rangle\oplus\Z\langle y\rangle$ 
which is defined by 
$$
g\cdot x=x,\;\; g\cdot y=-y+x.
$$
One can check easily that the integral $\Z_2$-representation on the
quotient modulo $\Z\langle x\rangle$, which is of cyclotomic type, 
does not lift to a summand of the same type in 
$\Z\langle x\rangle\oplus\Z\langle y\rangle$ because 
$\Z\langle x\rangle\oplus\Z\langle y\rangle=\Z[\Z_2]\langle(x-y)\rangle$ 
is of regular type. 

Likewise, a summand of cyclotomic type in a $\Z[\Z_p]$-submodule of a
$\Z[\Z_p]$-module may not be a summand of the same type in the 
$\Z[\Z_p]$-module.

\end{remark}

We end this section with a proof of Theorem 1.4. 
Recall that by Bryan \cite{Bry}, a smooth 
involution $g:X_\alpha\rightarrow X_\alpha$ is of odd type if and only if
$b_2^{+}(X_\alpha/\langle g\rangle)=1$. On the other hand, one can easily check that this
condition implies that one of the classes $[T_1]$, $[T_2]$ and
$[T_3]$ must be fixed by $g$. We assume without loss of generality that
$g^\ast[T_1]=[T_1]$. 

\begin{lemma}
Let $\Sigma$ be a non-spherical fixed component of $g$. Then (1)
$\chi(\Sigma)+\Sigma^2=0$, and (2) $\Sigma\cdot [T_1]=0$. 
\end{lemma}

\begin{proof}
First of all, let $\{\Sigma_j\}$ be the set of fixed components of $g$. Then
the Lefschetz fixed point theorem and the $G$-signature theorem (cf.
Theorem 3.4 and Theorem 3.6) imply that
$$
\left\{\begin{array}{ccc}
2+ t-(22-t) & = & \sum_j\chi(\Sigma_j)\\
2(2-t) & = & -16 + \sum_j \frac{2^2-1}{3}\cdot\Sigma_j^2,\\
\end{array} \right .
$$
where $t$ denotes the dimension of the $1$-eigenspace of $g$ in 
$H^2(X_\alpha;\R)$. It follows easily from the above equations that 
$\sum_j(\chi(\Sigma_j)+\Sigma_j^2)=0$. 

Now let $\{\Sigma_i\}$ be the set of components in $\{\Sigma_j\}$
such that $\Sigma_i^2<0$, and let $\{\Sigma_k\}$ be the set of 
components with $\Sigma_k^2\geq 0$. Then since $2d_1[T_1]$ is a Seiberg-Witten
basic class, by the generalized adjunction inequality, 
$$
\text{genus}(\Sigma_k)\geq 1+\frac{1}{2}(|2d_1[T_1]\cdot \Sigma_k|+\Sigma_k^2) 
$$
for each $k$. On the other hand, since $X_\alpha$ is even, $\Sigma_i^2
\leq -2$ for each $i$, so that 
$$
\text{genus}(\Sigma_i)\geq 1+\frac{1}{2}\Sigma_i^2. 
$$
Putting these two inequalities together, and with 
$\sum_j(\chi(\Sigma_j)+\Sigma_j^2)=0$, we obtain
$$
\sum_k |2d_1[T_1]\cdot \Sigma_k|\leq 0,
$$
which implies that $[T_1]\cdot\Sigma_k=0$ for each $k$ and 
$\chi(\Sigma_j)+\Sigma_j^2=0$ for each $j$. The lemma follows 
immediately. 

\end{proof}

\noindent{\bf Proof of Theorem 1.4}

\vspace{3mm}

Let $\Sigma$ be a fixed component of $g$ with genus $\geq 1$. We shall
prove that $\Sigma$ must be a torus of self-intersection $0$ and that
the class $[\Sigma]$ is a multiple of $[T_1]$ over $\Q$. Theorem 1.4 
follows easily from this and the result of Edmonds stated in
Proposition 3.2. 

We fix a $g$-equivariant decomposition $H^2(X_\alpha;\R)=H^{+}\oplus
H^{-}$ where $H^{+}$, $H^{-}$ are positive definite and negative definite
respectively. Since $b_2^{+}(X_\alpha/\langle g\rangle)=1$, there is a $1$-dimensional
subspace of $H^{+}$ which is fixed under $g$. We fix a vector $u\in H^{+}$
in this subspace such that $u^2=1$. Now because both $[T_1]$ and 
$[\Sigma]$ are fixed under $g$, we may write
$$
[T_1]=a_1u+\beta_1 \mbox{ and } [\Sigma]=a_2 u+\beta_2
$$
for some $a_1,a_2\in \R$ and $\beta_1,\beta_2\in H^{-}$. Note that 
our claim is trivially true if $[\Sigma]=0$. Assuming $[\Sigma]\neq 0$,
and note that $[T_1]\neq 0$, $T_1^2=0$, and $\Sigma^2=-\chi(\Sigma)\geq 0$, 
we must have $a_1,a_2\neq 0$. We may assume without loss of generality 
that $a_1,a_2>0$. With this understood, $T_1^2=0$, 
$\Sigma^2=-\chi(\Sigma)\geq 0$ and $[T_1]\cdot\Sigma=0$ give rise to
$$
a_1^2+\beta_1^2=0, \;\; a_2^2+\beta_2^2\geq 0, \mbox{ and }
a_1a_2+\beta_1\cdot\beta_2=0. 
$$
It follows easily that 
$$
|\beta_1\cdot\beta_2|=a_1a_2\geq (|\beta_1^2|\cdot |\beta_2^2|)^{1/2},
$$
which implies by the triangle inequality 
that $\beta_1,\beta_2$ must be linearly dependent and that
the above must hold with equality. 
It follows easily that $[\Sigma]$ is a multiple of $[T_1]$, and that
$\Sigma$ is a torus of self-intersection $0$.

\hfill $\Box$

\section{Proof of Theorem 1.7 --- Involutions in $\text{Aut}(E_8)$}

The proof of Theorem 1.7 requires a digression on the conjugacy classes 
of elements of order $2$ in $\text{Aut}(E_8)$. We shall give a brief
review of this material next, which is taken from Carter \cite{Car}.

Let $e_1,e_2,\cdots,e_8$ be a standard basis of $\R^8$, i.e., 
$(e_i,e_j)=\delta_{ij}$. Then the $E_8$ lattice is the lattice generated 
by the set of vectors 
$$
\Phi\equiv \{\pm e_i\pm e_j,\frac{1}{2}\sum_{i=1}^8 \epsilon_i e_i|
\epsilon_i=\pm 1, \prod_{i=1}^8\epsilon_i=1\}.
$$
Furthermore, $\Phi$ forms the root system of $E_8$. 

For any root $r\in\Phi$, there is an associated reflection 
$w_r\in \text{Aut}(E_8)$ defined by 
$$
w_r(x)=x-(r,x)r. 
$$
It is known that $\text{Aut}(E_8)$ is generated by $\{w_r|r\in\Phi\}$.

According to Lemma 5 of \cite{Car}, every involution $v\in \text{Aut}(E_8)$
can be written as a product $v=w_{r_1}\cdot w_{r_2}\cdot\cdots\cdot w_{r_k}$, 
where $k=l(v)$ equals the number
of $-1$-eigenvectors of $v$ in $\R^8$, and $r_1,\cdots,r_k$ are mutually
orthogonal roots. In particular, by changing $v$ to $(-1)\cdot v$ if 
necessary, we may assume that $k=l(v)\leq 4$. 

Let $f_1=e_1-e_2$, $f_2=e_2-e_3$, $\cdots$, $f_6=e_6-e_7$, $f_7=e_7+e_8$,
and 
$$
f_8=\frac{1}{2}(-e_1-e_2-e_3-e_4-e_5+e_6+e_7-e_8).
$$
Then one can easily check that $f_1,f_2,\cdots, f_8$ form a standard basis
for the $E_8$ lattice. In particular, $f_1,f_3,f_5,f_7$ are mutually 
orthogonal roots. 

Now according to \cite{Car} (see Lemma 11, Lemma 27, and Corollary (iv) 
following Proposition 38 in \cite{Car}), an involution $v\in \text{Aut}(E_8)$
is conjugate to one of 
$$
w_{f_1}, \;\;\; w_{f_1}\cdot w_{f_3}, \;\;\; w_{f_1}\cdot w_{f_3}\cdot w_{f_5}
$$
if $l(v)\leq 3$, and when $l(v)=4$, $v$ has two different conjugacy 
classes represented by 
$$
w_{f_1}\cdot w_{f_3}\cdot w_{f_5}\cdot w_{f_7} \mbox{ and }
w_{f_1}\cdot w_{f_3}\cdot w_{f_5}\cdot w_{f_7^\prime},
$$
where $f_7^\prime=e_7-e_8$. End of digression.

The following lemma is the starting point of our analysis.

\begin{lemma}
Suppose $\tau:X_\alpha\rightarrow X_\alpha$ is a smooth involution which
fixes the classes $[T_1]$, $[T_2]$ and $[T_3]$ and such that $\Theta(\tau)
=(v_1,v_2)\in \text{Aut}(E_8)\times \text{Aut}(E_8)$. {\em(}Recall that $\Theta$ 
is the homomorphism in Lemma 4.3.{\em)} Then both $v_1,v_2$ are conjugate to
the involution $w_{f_1}\cdot w_{f_3}\cdot w_{f_5}\cdot w_{f_7^\prime}$
in $\text{Aut}(E_8)$.  
\end{lemma}

\begin{proof}
Since $\tau$ fixes the classes $[T_1]$, $[T_2]$ and $[T_3]$, 
$b_2^{+}(X_\alpha/\langle\tau\rangle)=3$ by Lemma 4.1. Consequently $\tau$ is an
even type involution with $8$ isolated fixed points by Bryan \cite{Bry}.

The key property of $\tau$ we need here is that for any 
$x\in H_2(X_\alpha;\Z)$, the intersection product of $x$ with
$\tau\cdot x$ is even (cf. \cite{E1}). To see this, represent $x$
by a smooth surface $\Gamma$ in $X_\alpha$ which is away from the
fixed-point set of $\tau$, then perturb $\Gamma$ slightly so that
$\Gamma$ and $\tau(\Gamma)$ intersect transversely. 
It is easily seen that the intersection points of $\Gamma$ 
and $\tau(\Gamma)$ are paired up by $\tau$, and hence the claim. 
 
With this understood, we observe that if $v\in \text{Aut}(E_8)$ is an 
involution which is conjugate to any of the following $4$ involutions
$$
w_{f_1}, \;\;\; w_{f_1}\cdot w_{f_3},\;\;\; w_{f_1}\cdot w_{f_3}\cdot w_{f_5},
\;\;\; w_{f_1}\cdot w_{f_3}\cdot w_{f_5}\cdot w_{f_7},
$$
then there exists a root $x\in\Phi$ such that $(v(x),x)=1$. 
It suffices to check this for the above standard representatives of the 
conjugacy classes, which is done below.

If $v=w_{f_1}$, we take $x=f_2$. Then 
$(v(x),x)=(f_2+f_1,f_2)=2-1=1$. If $v=w_{f_1}\cdot w_{f_3}$, we take
$x=f_4$. Then $(v(x),x)=(f_4+f_3,f_4)=2-1=1$.
If $v=w_{f_1}\cdot w_{f_3}\cdot w_{f_5}$, we take $x=f_6$, and 
$(v(x),x)=(f_6+f_5,f_6)=1$. If 
$v=w_{f_1}\cdot w_{f_3}\cdot w_{f_5}\cdot w_{f_7}$, we take $x=f_8$,
and $(v(x),x)=(f_8+f_5,f_8)=1$. (On the other hand, 
if $v=w_{f_1}\cdot w_{f_3}\cdot w_{f_5}\cdot w_{f_7^\prime}$, then
a direct check shows that $(v(f_i),f_i)=0\pmod{2}$ for any 
$1\leq i\leq 8$.)

Consequently, by the classification of conjugacy classes of involutions in 
$\text{Aut}(E_8)$, we conclude that the involutions $v_1,v_2$ in 
$\Theta(\tau)=(v_1,v_2)\in\text{Aut}(E_8)\times \text{Aut}(E_8)$
have the following possibilities: either conjugate to
$w_{f_1}\cdot w_{f_3}\cdot w_{f_5}\cdot w_{f_7^\prime}$ or equal to 
$1$, $-1$.

It remains to show that neither of $v_1$, $v_2$ can be $1$ or $-1$. 
To this end, recall that there are classes $w_i\in H_2(X_\alpha;\Z)$,
$i=1,2,3$, which are dual to $[T_i]$. Moreover, since $\tau$ fixes $[T_1]$,
$[T_2]$ and $[T_3]$, for each $i=1,2,3$, there are 
$u_i\in \text{Span}([T_1],[T_2],[T_3])$ and 
$\alpha_i\in\text{Span}(\xi_k,\eta_k|1\leq k\leq 8)$ 
(the classes $\xi_k$, $\eta_k$ are defined in Lemma 4.3) such that
$$
\tau\cdot w_i=w_i+u_i+\alpha_i.
$$
(See the proof of Lemma 4.3 for details.) It follows easily from this that
$$
tr(\tau)|_{H_2(X_\alpha;\Z)}=6+tr(v_1)+tr(v_2).
$$
On the other hand, $\tau$ has $8$ isolated fixed points, so by the Lefschetz
fixed point theorem (cf. Theorem 3.4), $tr(\tau)|_{H_2(X_\alpha;\Z)}=8-2=6$, 
which implies that
$$
tr(v_1)+tr(v_2)=0. 
$$
Consequently, if $v_1=1$ or $-1$, $v_2$ must be $-1$ or $1$ respectively. 
Without loss of generality, we assume that $v_1=1$ and $v_2=-1$. 

Now $v_1=1$ means that $\tau$ acts trivially on 
$$
\text{Span}(\xi_k,[T_j]|1\leq k\leq 8, 1\leq j\leq 3)
/\text{Span}([T_1],[T_2],[T_3]).
$$
As we argued in the proof of Lemma 4.3, this implies that each $\xi_k$ is
fixed under $\tau$. We thus obtain a $\tau$-invariant decomposition 
$$
H_2(X_\alpha;\Z)=\text{Span}(\xi_k|1\leq k\leq 8)\oplus 
\text{Span}(\xi_k|1\leq k\leq 8)^\perp.
$$
(Here we use the fact that the intersection form on 
$\text{Span}(\xi_k|1\leq k\leq 8)$ is $-E_8$ which is unimodular.)
Suppose $\text{Span}(\xi_k|1\leq k\leq 8)^\perp=\Z[\Z_2]^r\oplus\Z^t\oplus
\Z[\mu_2]^s$ is a decomposition of the $\Z_2$-integral representation into
a block sum of summands of regular, trivial and cyclotomic types. Then 
correspondingly we have a decomposition 
$$
H_2(X_\alpha;\Z)=\Z[\Z_2]^r\oplus\Z^{t+8}\oplus\Z[\mu_2]^s.
$$
Now $tr(\tau)|_{H_2(X_\alpha;\Z)}=6$ implies $t+8-s=6$, which implies 
that $s=t+2>0$. However, $\tau$ is pseudofree and has a nonempty fixed-point
set, so that $s$ must be $0$ by Edmonds' result (cf. Prop. 3.1). This is a
contradiction, and the lemma follows. 

\end{proof}

Lemma 5.1 suggests that one should study $2$-subgroups of $\text{Aut}(E_8)$
whose elements of order $2$ are conjugate to 
$w_{f_1}\cdot w_{f_3}\cdot w_{f_5}\cdot w_{f_7^\prime}$. To this end we
need to recall a natual subgroup of $\text{Aut}(E_8)$ which contains all
the $2$-subgroups up to conjugacy. 

Consider the following two subgroups $H_0$ and $H_1$ 
of $\text{Aut}(E_8)$, 
$$
H_0=\{(\epsilon_i)|1\leq i\leq 8, \epsilon_i=\pm 1, 
\prod_{i=1}^8\epsilon_i=1\}\cong (\Z_2)^7,
$$
where $H_0$ acts by coordinate-wise multiplications on $\R^8$ with
respect to the standard basis $e_1,\cdots,e_8$, and 
$$
H_1=\{\sigma|\sigma \mbox{ is a permutation of } e_1,\cdots,e_8\}\cong S_8.
$$
Let $H\subset \text{Aut}(E_8)$ be the subgroup generated by $H_0$ and $H_1$.
Then $H$ is a semi-direct product of $H_0$ by $H_1$ with relations
$$
(\epsilon_i)\sigma=\sigma (\epsilon_i^\prime), \mbox{ where } 
\epsilon_i^\prime=\epsilon_{\sigma(i)},\;\;  \forall \;
(\epsilon_i)\in H_0,\;\;\sigma\in H_1.
$$
In particular, the Sylow $2$-subgroups of $H$ has order $2^{14}$ which is
the same as the order of Sylow $2$-subgroups of $\text{Aut}(E_8)$. Thus
by Sylow's theorem, up to conjugacy in $\text{Aut}(E_8)$ any $2$-subgroup
of $\text{Aut}(E_8)$ is contained in $H$. 

\begin{lemma}
(1) Let $v=(\epsilon_i)\hat{v}\in H$ where 
$(\epsilon_i)\in H_0$ and $\hat{v}\in H_1$. Then $v$ is conjugate to
$w_{f_1}\cdot w_{f_3}\cdot w_{f_5}\cdot w_{f_7^\prime}$ in $\text{Aut}(E_8)$
if and only if the following conditions are satisfied:
\begin{itemize}
\item either $\hat{v}=1$ or $\hat{v}=\sigma_1\sigma_2\sigma_3\sigma_4$ 
where $\sigma_i$ are disjoint transpositions,
\item $\epsilon_i=\epsilon_{\hat{v}(i)}$ for any $i$ and
$\#\{i|\epsilon_i=-1\}=0\pmod{4}$, 
\item when $\hat{v}=1$, $\#\{i|\epsilon_i=-1\}=4$. 
\end{itemize}

(2) Let $v=(\epsilon_i)\hat{v}\in H$ be of order $4$ such that $v^2$ is
conjugate to $w_{f_1}\cdot w_{f_3}\cdot w_{f_5}\cdot w_{f_7^\prime}$ in 
$\text{Aut}(E_8)$. Then there are the following two possibilities for $v$:
\begin{itemize}
\item Case (i) $\hat{v}^2=1$, where up to conjugacy in $H$, either
\begin{itemize}
\item $\hat{v}=(12)(34)$ with $\epsilon_1\epsilon_2
=\epsilon_3\epsilon_4=-1$, or 
\item $\hat{v}=(12)(34)(56)$ with  $\epsilon_1\epsilon_2
=\epsilon_3\epsilon_4=-1$, and $\epsilon_5\epsilon_6=1$, or
\item $\hat{v}=(12)(34)(56)(78)$ with $\epsilon_1\epsilon_2
=\epsilon_3\epsilon_4=-1$, and $\epsilon_5\epsilon_6=\epsilon_7\epsilon_8=1$.
\end{itemize}
\item Case (ii) $\hat{v}^2\neq 1$, where up to conjugacy in $H$,
$\hat{v}=\sigma_1\sigma_2$ for two disjoint $4$-cycles 
$\sigma_1,\sigma_2\in H_1$, and
$(\epsilon_i)$ satisfies $\epsilon_i=\epsilon_{\hat{v}(i)}$ for any $i$. 
\end{itemize}
\end{lemma}

\begin{proof}
(1) Note that up to conjugacy 
$w_{f_1}\cdot w_{f_3}\cdot w_{f_5}\cdot w_{f_7^\prime}$ may be 
characterized as the only involution $v\in \text{Aut}(E_8)$ such that
$v\neq -1$ and $(v(r),r)=0\pmod{2}$, $\forall r\in\Phi$. 

Now suppose $v=(\epsilon_i)\hat{v}\in H$ is an involution. Then
$$
1=v^2=(\epsilon_i)\hat{v}(\epsilon_i)\hat{v}=
(\epsilon_i)(\epsilon_{\hat{v}^{-1}(i)})\hat{v}\hat{v},
$$
which implies $\hat{v}^2=1$ and $\epsilon_i=\epsilon_{\hat{v}(i)}$ 
for any $i$.

Next we show that if $\hat{v}\neq 1$, then $\hat{v}$ must be a 
product of $4$ disjoint transpositions. To see this, suppose there 
exist $i\neq j$ such that $\hat{v}(i)=j$ and there exists a $k(\neq i,j)$ 
such that $\hat{v}(k)=k$, then for the root $e_i+e_k\in\Phi$,
$$
(v(e_i+e_k),(e_i+e_k))=(\epsilon_j e_j+\epsilon_k e_k, e_i+e_k)=\epsilon_k.
$$
Hence if $v$ is conjugate to 
$w_{f_1}\cdot w_{f_3}\cdot w_{f_5}\cdot w_{f_7^\prime}$, then either
$\hat{v}\neq 1$ or  $\hat{v}$ is a 
product of $4$ disjoint transpositions.

To see that $\#\{i|\epsilon_i=-1\}=0\pmod{4}$, note that
\begin{eqnarray*}
(v(f_8),f_8) & = & \frac{1}{4}(\#\{i|\epsilon_i=1\}-\#\{i|\epsilon_i=-1\})\\
             & = & \frac{1}{4}(8-2\#\{i|\epsilon_i=-1\})\\
             & = & 2-\frac{1}{2}\#\{i|\epsilon_i=-1\}. 
\end{eqnarray*}
Thus $(v(f_8),f_8)=0\pmod{2}$ if and only if 
$\#\{i|\epsilon_i=-1\}=0\pmod{4}$. When $\hat{v}=1$, $v=(\epsilon_i)$.
Since $v\neq 1$ or $-1$, we must have $\#\{i|\epsilon_i=-1\}=4$. 

Now suppose $v=(\epsilon_i)\hat{v}$ which satisfies the conditions in
(1) of the lemma. If $\hat{v}=1$, then $v=(\epsilon_i)$ where
$\#\{i|\epsilon_i=-1\}=4$. In particular, $v\neq 1, -1$. Moreover, 
for any root $r=\pm e_i\pm e_j$,
$$
(v(r),r)=\epsilon_i +\epsilon_j=0\pmod{2},
$$
and for any root $r=\frac{1}{2}(\sum_i \pm e_i)$,
$$
(v(r),r)=\frac{1}{4}(\#\{i|\epsilon_i=1\}-\#\{i|\epsilon_i=-1\})=0. 
$$
Hence $v$ is conjugate to 
$w_{f_1}\cdot w_{f_3}\cdot w_{f_5}\cdot w_{f_7^\prime}$ by the 
characterization of $w_{f_1}\cdot w_{f_3}\cdot w_{f_5}\cdot w_{f_7^\prime}$
we mentioned at the beginning of the proof.
When $\hat{v}=\sigma_1\sigma_2\sigma_3\sigma_4$ 
where $\sigma_i$ are disjoint transpositions, the conditions 
$\epsilon_i=\epsilon_{\hat{v}(i)}$ for any $i$ and
$\#\{i|\epsilon_i=-1\}=0\pmod{4}$ imply that $v$ is conjugate to
$\hat{v}=\sigma_1\sigma_2\sigma_3\sigma_4$ by an element of $H_0$.
On the other hand, $\hat{v}=\sigma_1\sigma_2\sigma_3\sigma_4$ is clearly
conjugate to $(12)(34)(56)(78)$ in $H_1$, which is exactly 
$w_{f_1}\cdot w_{f_3}\cdot w_{f_5}\cdot w_{f_7^\prime}$. This proves 
part (1) of the lemma.

(2) Let $v=(\epsilon_i)\hat{v}$ be of order $4$ such that $v^2$ is 
conjugate to $w_{f_1}\cdot w_{f_3}\cdot w_{f_5}\cdot w_{f_7^\prime}$.
We have $v^2=(\epsilon_i)(\epsilon_{\hat{v}^{-1}(i)})\hat{v}^2$.

Case (i) where $\hat{v}^2=1$. Then by part (1) and up to conjugacy by
an element of $H_1$,  
$$
v^2=(\epsilon_i)(\epsilon_{\hat{v}^{-1}(i)})\hat{v}^2=(-1,-1,-1,-1,1,1,1,1). 
$$
This implies that for $i=1,2,3,4$, $\epsilon_i\epsilon_{\hat{v}(i)}=-1$,
and in particular, $\hat{v}(i)\neq i$. 
Up to further conjugation by an element of $H_1$, $\hat{v}$ has the 
following three possibilities: 
$$
(12)(34), (12)(34)(56), (12)(34)(56)(78).
$$
The corresponding conditions that $(\epsilon_i)$ must satisfy follow directly
from 
$$
v^2=(\epsilon_i)(\epsilon_{\hat{v}^{-1}(i)})=(-1,-1,-1,-1,1,1,1,1).
$$

Case (ii) where $\hat{v}^2\neq 1$. Then 
$v^2=(\epsilon_i)(\epsilon_{\hat{v}^{-1}(i)})\hat{v}^2$, which, as we have 
seen in part (1), is conjugate to a product of $4$ disjoint transpositions.
This implies that $\hat{v}$ is a product of $2$ disjoint $4$-cycles and
$(\epsilon_i)(\epsilon_{\hat{v}^{-1}(i)})=(1)$ implies that 
$\epsilon_i=\epsilon_{\hat{v}(i)}$ for any $i$. 

\end{proof}

We remark that with Lemma 5.2 one can easily show that any $2$-group
of order $\leq 8$ (including $Q_8$ in particular) as well as some other groups
of small order (e.g. $S_3$, $A_4$, or even $S_4$) can be realized as a 
subgroup of $H$ whose order $2$ elements are all conjugate to 
$w_{f_1}\cdot w_{f_3}\cdot w_{f_5}\cdot w_{f_7^\prime}$. (One may even
attempt to classify these subgroups of $H$.) With this 
understood, the following lemma provides the additional constraints 
needed for the case of $Q_8$ in Theorem 1.7.

\begin{lemma}
Suppose $Q_8$ acts on $X_\alpha$ smoothly, such that (1) the classes 
$[T_1]$, $[T_2]$ and $[T_3]$ are fixed under the action, (2) the actions 
by the elements of order $4$ of $Q_8$ are mutually conjugate. Then each
element of order $4$ of $Q_8$ has exactly $4$ isolated fixed points in 
$X_\alpha$.
\end{lemma}

\begin{proof}
Let $g\in Q_8$ be an order $4$ element. Since $[T_1]$, $[T_2]$ and $[T_3]$ 
are fixed under the action, $b_2^{+}(X_\alpha/\langle g\rangle)
=b_2^{+}(X_\alpha/\langle g^2\rangle)=3$,
and in particular, $g^2$ is an even type involution with $8$ isolated fixed
points. Since $\text{Fix}(g)\subset \text{Fix}(g^2)$, we see that $g$ has
at most $8$ isolated fixed points.

We shall prove next that the number of fixed points of $g$ is either $4$,
$6$ or $8$. To see this, note that the fixed points of $g$ fall into two
different classes according to their local representations. Denote by
$s_{+}$ the number of fixed points where the weights of the local 
representation are $(1,3)$ and denote by $s_{-}$ the number of fixed 
points where the weights are $(1,1)$ or $(3,3)$. Note that if $t$ 
is the dimension of the $1$-eigenspace of $g$ in $H^2(X_\alpha;\R)$, 
the dimension of the $(-1)$-eigenspace must be $14-t$, because $g^2$ 
has $8$ isolated fixed points so that the dimension of the 
$1$-eigenspace of $g^2$ is $14$. Now by the Lefschetz fixed point theorem 
(cf. Theorem 3.4) and the $G$-signature theorem (cf. \cite{HZ}), we have
$$
\left\{\begin{array}{ccc}
2+t-(14-t) & = & s_{+}+s_{-}\\
4(6-t) & = & -16+ 2s_{+}+(-2)s_{-}.\\
\end{array} \right .
$$
Here we use the fact that $b_2^{+}(X_\alpha/g)=3$, and the fact that 
the signature defect at a fixed point of $g$ of type $(1,3)$ and type $(1,1)$ 
or $(3,3)$ is $2,-2$ respectively, and the signature defect at a fixed 
point of $g^2$ is $0$. The solutions for $s_{+}$, $s_{-}$ (note that 
$s_{+}+s_{-}\leq 8$) are $s_{+}=4$ and $s_{-}=0,2$ or $4$. Our claim
about the number of fixed points of $g$ follows immediately.

Now $Q_8=\{i,j,k|i^2=j^2=k^2=-1,\; ij=k,\; jk=i,\; ki=j\}$ where by the
assumption the actions by $i,j,k$ are all conjugate. In particular,
they have the same number of fixed points, which is either $4$, $6$ or $8$.
Suppose $i,j,k$ all have $8$ fixed points. Then each of the $8$ fixed points
of $-1$ is also fixed by the entire group $Q_8$. But one of them must be
a fixed point of $i$ of type $(1,1)$ or $(3,3)$, which, however, 
contradicts the relation $j^{-1}ij=i^{-1}$. Hence $i,j,k$ can not have
$8$ fixed points each. Suppose $i,j,k$ each has $6$ fixed points. Then
$i,j$ each fixes $6$ of the $8$ fixed points of $-1$, so that
they must have $4$ common fixed points, which should also be fixed by
$k=ij$. Let $x_1,\cdots,x_4$ denote these $4$ common fixed points, and
let $x_5$,$x_6$ denote the other $2$ fixed points of $i$ which is not 
fixed by $j$, and let $x_7$, $x_8$ be the remaining $2$ points which 
are not fixed by $i$. Since $j^{-1}ij=i^{-1}$, $j$ has to switch $x_5$, 
$x_6$, so that $j$ must fix both $x_7$, $x_8$ because $j$ has $6$ fixed 
points. It follows easily that $k=ij$ does not fix any of the points 
$x_5$, $x_6$, $x_7$, $x_8$, which contradicts the assumption that
$k$ also has $6$ fixed points. Hence $i,j,k$ each has $4$ fixed points,
and the lemma follows. 

\end{proof}

\noindent{\bf Proof of Theorem 1.7}

\vspace{3mm}

Let $G$ be a finite group acting smoothly and effectively on $X_\alpha$. Then the classes
$[T_1]$, $[T_2]$, $[T_3]$ are fixed by the commutator subgroup $[G,G]$.
Moreover, under the homomorphism $\Theta:G\rightarrow \text{Aut}(E_8\oplus
E_8)$, $[G,G]$ is mapped into the subgroup $\text{Aut}(E_8)\times 
\text{Aut}(E_8)$ of index $2$. We denote by $\Theta_i$, $i=1,2$,  the  
homomorphisms into $\text{Aut}(E_8)$ such that $\Theta=(\Theta_1,\Theta_2)$
on $[G,G]$.

Now suppose $K$ is a subgroup of $[G,G]$ which is isomorphic to $(\Z_2)^4$.
Then by Lemma 5.1, for every element $g\in K$ such that $g\neq 1$,
$\Theta_1(g)\in\text{Aut}(E_8)$ is conjugate to $w_{f_1}\cdot w_{f_3}
\cdot w_{f_5}\cdot w_{f_7^\prime}$. In particular, the trace of 
$\Theta_1(g)$ equals $0$ (cf. Lemma 5.2 (1)).  Now consider the 
$8$-dimensional representation $V$ of $K$ induced by $\Theta_1$. We have
$$
\dim V^K=\frac{1}{|K|}\sum_{g\in K}tr(g)
=\frac{1}{|K|}(8+\sum_{1\neq g\in K}tr(g))=\frac{8}{16},
$$ 
which is a contradiction. This proves that if $[G,G]$ contains $(\Z_2)^4$
as a subgroup, then $G$ can not act smoothly and effectively on $X_\alpha$. 

Next we assume $K\subset [G,G]$ is a subgroup isomorphic to $Q_8$, where
$$
Q_8=\{i,j,k|i^2=j^2=k^2=-1, ij=k,\;jk=i,\;ki=j\}.
$$
By the assumption in Theorem 1.7, the actions of order $4$ elements are
all conjugate, so that by Lemma 5.3, each of $i,j,k$ has $4$ isolated 
fixed points. It then follows easily from the Lefschetz fixed point 
theorem (cf. Theorem 3.4) that 
$$
tr(\Theta_1(g))+tr(\Theta_2(g))=-4, \;\;\mbox { where } g=i,j,k \in Q_8.
$$
On the other hand, each of $\Theta_1(g)$, $\Theta_2(g)$, $g=i,j,k$,
is an order $4$ element of $\text{Aut}(E_8)$ whose square is conjugate 
to $w_{f_1}\cdot w_{f_3}\cdot w_{f_5}\cdot w_{f_7^\prime}$.  
From the description in Lemma 5.2 (2), the trace of each of 
$\Theta_1(g),\Theta_2(g)$,
$g=i,j,k$, is even and is bounded between $-4$ and $4$. It follows that
for $l=1,2$, $tr(\Theta_l(g))=0,-2,-4$ where $g=i,j,k\in Q_8$. 

There are two possibilities which we will discuss separately. First,
suppose for $g$ equaling one of $i,j,k\in Q_8$, $tr(\Theta_1(g))=0$
or $-4$. Note that correspondingly $tr(\Theta_2(g))=-4$ or $0$. So
without loss of generality we may assume that $tr(\Theta_1(i))=-4$
(i.e. $g=i$). Then from the description in Lemma 5.2 (2), we must have
(up to conjugacy) 
$\Theta_1(i)=(\epsilon_i)\hat{v}$ with $\hat{v}=(12)(34)$,
$\epsilon_1\epsilon_2=\epsilon_3\epsilon_4=-1$, $\epsilon_l=-1$ for
$l=5,6,7,8$. Note also in this case we have 
$\Theta_1(-1)=(-1,-1,-1,-1,1,1,1,1)$. 

We discuss the possibilities 
for $tr(\Theta_1(j))$. Suppose $tr(\Theta_1(j))=-4$. Then 
$\Theta_1(j)=(\epsilon_i)\hat{v}$ where $\hat{v}$ is a product of $2$
disjoint transpositions with each of $5,6,7,8$ being fixed, and where 
$\epsilon_l=-1$ for $l=5,6,7,8$. It follows easily that 
$\Theta_1(k)=(\epsilon_i)\hat{v}$ where
$\hat{v}$ fixes $5,6,7,8$ and $\epsilon_l=1$ for $l=5,6,7,8$. But this implies 
that $tr(\Theta_1(k))=4$ which is a contradiction. Suppose 
$tr(\Theta_1(j))=-2$. Then up to conjugacy without effecting $\Theta_1(i)$,
$\Theta_1(j)=(\epsilon_i)\hat{v}$ where 
$\hat{v}=\sigma_1\sigma_2(56)(7)(8)$ and $\epsilon_7=\epsilon_8=-1$. 
It follows that $\Theta_1(k)=(\epsilon_i)\hat{v}$ with 
$\hat{v}=\sigma_1^\prime\sigma_2^\prime(56)(7)(8)$ but 
$\epsilon_7=\epsilon_8=1$. Consequently $tr(\Theta_1(k))=2$ 
which is also a contradiction. Finally, suppose $tr(\Theta_1(j))=0$. Then
either $\Theta_1(j)=(\epsilon_i)\hat{v}$ with $\hat{v}$ a product of 
$4$ disjoint transpositions or with $\hat{v}$ fixing each of $5,6,7,8$
and exactly two of $\epsilon_5, \epsilon_6,\epsilon_7,\epsilon_8$ equal
to $-1$. In any event, it follows that $tr(\Theta_1(k))=0$ also. 
Now $tr(\Theta_1(j))=tr(\Theta_1(k))=0$ implies that
$tr(\Theta_2(j))=tr(\Theta_2(k))=-4$, which is a case that has been already
shown impossible. Hence we have eliminated the first possibility that 
for $g$ equaling one of $i,j,k\in Q_8$, $tr(\Theta_1(g))=0$ or $-4$.

Next we consider the second possibility that for any $g=i,j,k\in Q_8$,
$tr(\Theta_1(g))=-2$. Then from the description in Lemma 5.2 (2), we see
that each $\Theta_1(g)=(\epsilon_i)\hat{v}$, where $\hat{v}$ is a product
of $3$ disjoint transpositions. This turns out to be impossible,
because if we assume without loss of generality that $\Theta_1(i)
=(\epsilon_i)\hat{v}$ where $\hat{v}=(12)(34)(56)$
(and $\Theta_1(-1)=(-1,-1,-1,-1,1,1,1,1)$), then $\Theta_1(j)
=(\epsilon_i)\hat{v}$ with $\hat{v}=\sigma_1\sigma_2(56)(7)(8)$ or
$\hat{v}=\sigma_1\sigma_2(5)(6)(78)$, which implies that 
$\Theta_1(k)=(\epsilon_i)\hat{v}$ where 
$\hat{v}=\sigma_1^\prime\sigma_2^\prime(5)(6)(7)(8)$ or 
$\hat{v}=\sigma_1^\prime\sigma_2^\prime(56)(78)$ respectively.
In any event, it contradicts the assumption that $tr(\Theta_1(k))=-2$,
and hence our claim. This proves that if $[G,G]$ 
contains $Q_8$ as a subgroup, then $G$ can not act smoothly and
effectively on $X_\alpha$

The proof of Theorem 1.7 is thus completed. 

\hfill $\Box$

\begin{corollary}
The following maximal symplectic $K3$ groups 
$$
G=M_{20}, F_{384}, A_{4,4}, T_{192}, H_{192}, T_{48}
$$
can not act smoothly and effectively on $X_{\alpha}$. 
\end{corollary}

\begin{proof}
The structures of these groups and their commutator subgroups are listed 
below (cf. \cite{Mu, Xiao}):

\begin{itemize}
\item $G=M_{20}=2^4 A_5$, $[G,G]=G=2^4 A_5$.
\item $G=F_{384}=4^2 S_4$, $[G,G]=4^2 A_4$.
\item $G=A_{4,4}=2^4 A_{3,3}$, $[G,G]=A_4^2$.
\item $G=T_{192}=(Q_8\ast Q_8)\times_\phi S_3$,
$[G,G]=(Q_8\ast Q_8)\times_\phi \Z_3$.
\item $G=H_{192}=2^4 D_{12}$, $[G,G]=2^4\Z_3$.
\item $G=T_{48}=Q_8\times_\phi S_3$, $[G,G]=T_{24}=Q_8\times_\phi \Z_3$.
\end{itemize}
The corollary is evident for all cases except for the case where
$G=F_{384}$. We shall prove that in this case $[G,G]=4^2 A_4$ contains 
a subgroup isomorphic to $(\Z_2)^4$. To this end, we recall the structure 
of $F_{384}$ 
(cf. \cite{Mu}, pages 190-191). $F_{384}=4^2 S_4$ is a semi-direct product
of $(\Z_4)^2$ by $S_4$, where 
$(\Z_4)^2=\{(a,b,c,d)|a+b+c+d=0\}\subset (\Z_4)^4$ modulo the diagonal
subgroup. The action of $S_4$ is given by permutations of the $4$ coordinates.
One can check directly that the $(\Z_2)^2$-subgroup of $(\Z_4)^2$ generated
by $(2,2,0,0)$ and $(2,0,2,0)$ are fixed under the action of 
$(12)(34), (13)(24)\in A_4\subset S_4$, hence the commutator 
$[G,G]=4^2 A_4$ contains a subgroup isomorphic to $(\Z_2)^4$. 

\end{proof}

\section{Symplectic cyclic actions}

In this section we prove Theorem 1.8. The proof draws heavily on 
our previous work \cite{CK} concerning the fixed-point set structure
of a symplectic $\Z_p$-action on a minimal symplectic $4$-manifold with
$c_1^2=0$, which we shall recall first.

Let $\omega$ be an orientation compatible symplectic structure on $X_\alpha$, 
and let $G$ be a finite group acting on $X_\alpha$ which preserves $\omega$. 
Then by Lemma 2.3 (2), $G$ fixes the classes $[T_1], [T_2], [T_3]$, and 
therefore by Lemma 4.1, $G$ acts trivially on a $3$-dimensional subspace
of $H^{2}(X_\alpha;\R)$ which consists of elements of positive square. 
As we argued in \cite{CK}, a $G$-equivariant version 
of Taubes' work in \cite{T1, T2} applies here, so that for any $G$-equivariant 
$\omega$-compatible almost complex structure $J$, the canonical class
$c_1(K)$ is represented by a finite set of $J$-holomorphic curves 
$\{C_i\}$ with positive weights $\{n_i\}$, i.e., $c_1(K)=\sum_i n_i C_i$,
which has the following properties:
\begin{itemize}
\item The set $\cup_i C_i$ is $G$-invariant.
\item Any fixed point of $G$ in the complement of $\cup_i C_i$ is isolated 
with local representation contained in $SL_2(\C)$.
\end{itemize}

One may further analyze the rest of the fixed points through the induced 
action in a neighborhood of $\cup_i C_i$. To this end, it is useful to take 
note that the connected components of $\cup_i C_i$ may be divided 
into the following three types:
\begin{itemize}
\item [{(A)}] A single $J$-holomorphic curve of self-intersection $0$ 
which is either an embedded torus, or a cusp sphere, or a nodal sphere.
\item [{(B)}] A union of two embedded $(-2)$-spheres intersecting at a 
single point with tangency of order $2$.
\item [{(C)}] A union of embedded $(-2)$-spheres intersecting transversely. 
\end{itemize}

A type (C) component may be conveniently represented by one of the graphs 
of type $\tilde{A}_n$, $\tilde{D}_n$, $\tilde{E}_6$, $\tilde{E}_7$
or $\tilde{E}_8$ listed in Figure 2, where a vertex in a graph represents 
a $(-2)$-sphere and an edge connecting two vertices represents a transverse, 
positive intersection point of the two $(-2)$-spheres represented by the
vertices. 

\begin{figure}[ht]

\setlength{\unitlength}{1pt}

\begin{picture}(420,420)(-210,-195)

\linethickness{0.5pt}

\put(-170,165){$\tilde{A}_n$}

\put(-120,165){\circle*{3}}

\put(-120,165){\line(1,1){30}}

\put(-90,195){\circle*{3}}

\put(-90,195){\line(1,0){30}}

\put(-60,195){\circle*{3}}

\put(-60,195){\line(1,-1){30}}

\put(-30,165){\circle*{3}}

\put(-120,165){\line(1,-1){30}}

\put(-90,135){\circle*{3}}

\put(-90,135){\line(1,0){30}}

\put(-60,135){\circle*{3}}

\multiput(-60,135)(3,3){10}{\line(1,1){2}}

\put(20,165){($n+1$ vertices, $n\geq 1$)}

\put(-170,90){$\tilde{D}_n$}

\put(-120,105){\circle*{3}}

\put(-120,105){\line(2,-1){30}}

\put(-120,75){\circle*{3}}

\put(-120,75){\line(2,1){30}}

\put(-90,90){\circle*{3}}

\put(-90,90){\line(1,0){30}}

\put(-60,90){\circle*{3}}

\put(-30,90){\circle*{3}}

\multiput(-60,90)(3,0){10}{\line(1,0){2}}

\put(0,90){\circle*{3}}

\put(-30,90){\line(1,0){30}}

\put(0,90){\line(2,1){30}}

\put(0,90){\line(2,-1){30}}

\put(30,105){\circle*{3}}

\put(30,75){\circle*{3}}

\put(60,90){($n+1$ vertices, $n\geq 4$)}

\put(-170,30){$\tilde{E}_6$}

\put(-120,30){\circle*{3}}

\put(-120,30){\line(1,0){30}}

\put(-90,30){\circle*{3}}

\put(-90,30){\line(1,0){30}}

\put(-60,30){\circle*{3}}

\put(-60,30){\line(1,0){30}}

\put(-30,30){\circle*{3}}

\put(-30,30){\line(1,0){30}}

\put(0,30){\circle*{3}}

\put(-60,30){\line(0,-1){60}}

\put(-60,0){\circle*{3}}

\put(-60,-30){\circle*{3}}

\put(-170,-60){$\tilde{E}_7$}

\put(-120,-60){\circle*{3}}

\put(-120,-60){\line(1,0){180}}

\put(-90,-60){\circle*{3}}

\put(-60,-60){\circle*{3}}

\put(-30,-60){\circle*{3}}

\put(0,-60){\circle*{3}}

\put(30,-60){\circle*{3}}

\put(60,-60){\circle*{3}}

\put(-30,-60){\line(0,-1){30}}

\put(-30,-90){\circle*{3}}

\put(-170,-150){$\tilde{E}_8$}

\put(-120,-150){\circle*{3}}

\put(-120,-150){\line(1,0){210}}

\put(-90,-150){\circle*{3}}

\put(-60,-150){\circle*{3}}

\put(-30,-150){\circle*{3}}

\put(0,-150){\circle*{3}}

\put(30,-150){\circle*{3}}

\put(60,-150){\circle*{3}}

\put(90,-150){\circle*{3}}

\put(-60,-150){\line(0,-1){30}}

\put(-60,-180){\circle*{3}}

\end{picture}

\caption{}

\end{figure}

With the preceding understood, the following lemma is specially tailored 
for the present situation in order to control the number of type (B) or 
type (C) components.

\begin{lemma}
Let $J$ be any $\omega$-compatible almost complex structure on $X_\alpha$,
and let $c_1(K)=\sum n_iC_i$ where $\{C_i\}$ is a finite set of 
$J$-holomorphic curves and $n_i\geq 1$. Then each $C_i$ lies in the 
orthogonal complement of $\mbox{Span }([T_1],[T_2],[T_3])$.
\end{lemma}

\begin{proof}
The key point here is that some multiple of each $[T_j]$ can be
represented by $J$-holomorphic curves. The details of the proof go
as follows.

First of all, by Lemma 2.2, we may assume without loss of generality that
$c_1(K)=2(d_1[T_1]+d_2[T_2]+d_3[T_3])$. Because the classes
$-2(d_2[T_2]+d_3[T_3])$, $-2(d_1[T_1]+d_3[T_3])$ and
$-2(d_1[T_1]+d_2[T_2])$ are Seiberg-Witten basic classes, by the main
theorem of Taubes in \cite{T2}, for a generic $\omega$-compatible
almost complex structure $J^\prime$, each $d_j[T_j]$ for $j=1,2,3$ is
Poincar\'{e} dual to $\sum_{k=1}^{N_j} m_{j,k} \Gamma_{j,k}$, where
$m_{j,k}\geq 1$ are integers and $\Gamma_{j,k}$ are (connected) 
embedded $J^\prime$-holomorphic curves which are disjoint for each fixed 
$j$. Moreover, since $X_\alpha$ is minimal and $J^\prime$ is generic, 
all $\Gamma_{j,k}$ have nonzero 
genus. We further notice that for each fixed $j$ the numbers $N_j$ and 
$m_{j,k}$ and the genus of each $\Gamma_{j,k}$ are bounded by a constant 
independent of the almost complex structure $J^\prime$. We take a sequense 
of generic $J^\prime$ converging in $C^\infty$ to the given $J$, and by 
passing to a subsequence we may assume that $N_j$, $m_{j,k}$ and the 
genus of $\Gamma_{j,k}$ are independent of $J^\prime$ throughout.

By Gromov compactness theorem (cf. e.g. \cite{McDS}), 
each $\Gamma_{j,k}$ converges to
a limit $\sum_{l=1}^{M_{j,k}}n_{j,k,l}C_{j,k,l}$ where
each $C_{j,k,l}$ is a (nonconstant) $J$-holomorphic curve,
$n_{j,k,l}\geq 1$ and $\cup_{l=1}^{M_{j,k}}C_{j,k,l}$ is connected.
Note that
$$
c_1(K)=2(d_1[T_1]+d_2[T_2]+d_3[T_3])=2\sum_{j=1}^3\sum_{k=1}^{N_j}
\sum_{l=1}^{M_{j,k}}m_{j,k} n_{j,k,l}C_{j,k,l}.
$$
Furthermore, the fact that $c_1(K)^2=0$ and $X_\alpha$ is minimal allows 
us to analyze the structure of
$\cup_{j=1}^3\cup_{k=1}^{N_j}\cup_{l=1}^{M_{j,k}}C_{j,k,l}$,
as shown in \cite{CK}. In particular, the connected
components of the union $\cup_{j=1}^3\cup_{k=1}^{N_j}\cup_{l=1}^{M_{j,k}}
C_{j,k,l}$ may be divided into the following three types (the 
classification differs slightly from the one we mentioned earlier):
\begin{itemize}
\item [{(a)}] A single $J$-holomorphic curve of self-intersection $0$.
\item [{(b)}] A union of two embedded $(-2)$-spheres.
\item [{(c)}] A union of at least three embedded $(-2)$-spheres intersecting
transversely.
\end{itemize}

With the above preparation, we shall prove next that each $C_i$ lies in the
orthogonal complement of $\mbox{Span }([T_1],[T_2],[T_3])$.
It suffices to show that for each $(j,k)$, $\Gamma_{j,k}\cdot C_i=0$, or
equivalently $(\sum_{l=1}^{M_{j,k}}n_{j,k,l}C_{j,k,l})\cdot C_i=0$.
We begin by recalling that $c_1(K)\cdot C_i=0$ (cf. \cite{CK}, Lemma 3.3),
so that $C_i$ is either disjoint from
$\cup_{j=1}^3\cup_{k=1}^{N_j}\cup_{l=1}^{M_{j,k}}C_{j,k,l}$, in which case
$(\sum_{l=1}^{M_{j,k}}n_{j,k,l}C_{j,k,l})\cdot C_i=0$ holds true
automatically, or $C_i$ is contained as one of the $J$-holomorphic 
curves $C_{j,k,l}$. 

At this point, we need to make use of the fact that $\Gamma_{j,k}^2=0$,
whose proof is postponed to the end of the proof of this lemma. Accepting 
it momentarily, we shall continue with the proof of the lemma. It is clear 
that we only need to verify the case where $\cup_{l=1}^{M_{j,k}}C_{j,k,l}$ 
and $C_i$ lie in the same component of 
$\cup_{j=1}^3\cup_{k=1}^{N_j}\cup_{l=1}^{M_{j,k}}C_{j,k,l}$.
Since each $\cup_{l=1}^{M_{j,k}}C_{j,k,l}$ is connected and
$(\sum_{l=1}^{M_{j,k}}n_{j,k,l}C_{j,k,l})^2=\Gamma_{j,k}^2=0$, 
it follows easily that if $C_i$ lies in a type (a) or (b) component of
$\cup_{j=1}^3\cup_{k=1}^{N_j}\cup_{l=1}^{M_{j,k}}C_{j,k,l}$, then
$(\sum_{l=1}^{M_{j,k}}n_{j,k,l}C_{j,k,l})\cdot C_i=0$ holds true.
It remains to check the case where $C_i$ lies in a type (c) component.
To this end we recall that a type (c) component corresponds to a graph
of type $\tilde{A}_n$, $\tilde{D}_n$, $\tilde{E}_6$, $\tilde{E}_7$
or $\tilde{E}_8$ as discussed in \cite{BPV}, Lemma 2.12 (ii). Each
graph defines a positive semi-definite quadratic form which is canonically
associated with the intersection form of the $J$-holomorphic curves 
$C_{j,k,l}$ in the type (c) component. The key property we will use here 
is that the positive semi-definite quadratic form has a $1$-dimensional 
annihilator. Now it is clear that 
$(\sum_{l=1}^{M_{j,k}}n_{j,k,l}C_{j,k,l})^2=0$ implies
that $\sum_{l=1}^{M_{j,k}}n_{j,k,l}C_{j,k,l}$ must be an annihilator for the
positive semi-definite quadratic form, which implies that
$(\sum_{l=1}^{M_{j,k}}n_{j,k,l}C_{j,k,l})\cdot C_i=0$.

We end the proof by showing that $\Gamma_{j,k}^2=0$. This follows from the
fact that $[T_j]^2=0$ by a standard argument involving the Sard-Smale
theorem and the adjunction formula for pseudoholomorphic curves (cf.
\cite{McDS}). The details are sketched below. The dimension of the moduli
space of $J^\prime$-holomorphic curves which contains $\Gamma_{j,k}$ equals
$d=2(-c_1(K)\cdot\Gamma_{j,k}+\text{genus}(\Gamma_{j,k})-1)$. (Here we use
the fact that $\Gamma_{j,k}$ has nonzero genus.) Since $\Gamma_{j,k}$ is
embedded, the adjunction formula 
$$
2\cdot\text{genus}(\Gamma_{j,k})-2=
\Gamma_{j,k}^2+c_1(K)\cdot\Gamma_{j,k}
$$ 
gives rise to
$d=\Gamma_{j,k}^2-c_1(K)\cdot\Gamma_{j,k}$. Now $J^\prime$ is chosen generic
so that $d\geq 0$ must hold, which implies that $\Gamma_{j,k}^2\geq 
c_1(K)\cdot\Gamma_{j,k}$. Again by the adjunction formula, we have
$$
\Gamma_{j,k}^2\geq\frac{1}{2}(\Gamma_{j,k}^2+c_1(K)\cdot\Gamma_{j,k})
=\text{genus}(\Gamma_{j,k})-1\geq 0.
$$ 
With this, $\Gamma_{j,k}^2=0$ follows easily from 
$(\sum_{k=1}^{N_j}m_{j,k}\Gamma_{j,k})^2=(d_j[T_j])^2=0$. 

\end{proof}

The preceding lemma has the following useful corollary. 
Let $\Lambda$ be a component of $\cup_i C_i$ of either type (B)
or type (C), and let $C$ be a $(-2)$-sphere in $\Lambda$. Recall that 
the orthogonal complement of $\mbox{Span }([T_1],[T_2],[T_3])$ is 
$$
\mbox{Span }([T_j],\xi_k,\eta_k|j=1,2,3, 1\leq k\leq 8)
$$ 
where $\xi_k,\eta_k$ are the classes in $H_2(X_\alpha;\Z)$ which correspond 
to the two standard bases of the $-E_8$ form defined in Lemma 4.2.
We denote by $\underline{C}$ the projection into
$$
\mbox{Span }([T_j],\xi_k,\eta_k|j=1,2,3, 1\leq k\leq 8)/
\mbox{Span }([T_1], [T_2], [T_3]).
$$
Since $C$ has nontrivial self-intersection, its projection
$\underline{C}$ must be nonzero. We denote by $L_\Lambda$ the sublattice
spanned by the projections of $(-2)$-spheres in $\Lambda$.

\begin{lemma}
For any component $\Lambda$, $L_\Lambda$ is contained in either
$$
\mbox{Span }([T_j],\xi_k|j=1,2,3, 1\leq k\leq 8)
/\mbox{Span }([T_1],[T_2],[T_3])
$$ 
or
$$
\mbox{Span }([T_j],\eta_k|j=1,2,3, 1\leq k\leq 8)
/\mbox{Span }([T_1],[T_2],[T_3]),
$$ 
and for any two distinct components $\Lambda,\Lambda^\prime$, the 
corresponding sublattices $L_\Lambda$, $L_{\Lambda^\prime}$ are 
orthogonal to each other. Moreover, if $\Lambda$ is of type (B),
then $L_\Lambda=\langle -2\rangle $ {\em(}i.e. $L_\Lambda$ is
a $A_1$-root lattice{\em)}, and if $\Lambda$ is of type (C), represented 
by a graph of type $\tilde{A}_n$, $\tilde{D}_n$, $\tilde{E}_6$, $\tilde{E}_7$
or $\tilde{E}_8$ listed in Figure 2, then $L_\Lambda$ is a root lattice of
the corresponding type {\em(}i.e. of type ${A}_n$, ${D}_n$, ${E}_6$, ${E}_7$
or ${E}_8${\em)}.
\end{lemma}

\begin{proof}
Let $C$ be a $(-2)$-sphere in $\Lambda$. Write $\underline{C}=\xi+\eta$,
where 
$$
\xi\in\mbox{Span }([T_j],\xi_k|j=1,2,3, 
1\leq k\leq 8)/\mbox{Span }([T_1],[T_2],[T_3])
$$
and 
$$
\eta\in\mbox{Span }([T_j],\eta_k|j=1,2,3, 1\leq k\leq 8)
/\mbox{Span }([T_1],[T_2],[T_3]).
$$
We claim that either $\xi$ or $\eta$ is zero. Suppose to the contrary that
neither of them is zero. Then since $-E_8$ is negative definite and even, 
$\xi^2, \eta^2\leq -2$, which implies that 
$\underline{C}^2=\xi^2+\eta^2\leq -4$. 
But this contradicts $\underline{C}^2=C^2=-2$, and the claim follows.

Now for each $(-2)$-sphere $C$ in $\Lambda$, its projection $\underline{C}$
lies in either 
$$
\mbox{Span }([T_j],\xi_k|j=1,2,3, 
1\leq k\leq 8)/\mbox{Span }([T_1],[T_2],[T_3])
$$
or 
$$
\mbox{Span }([T_j],\eta_k|j=1,2,3, 1\leq k\leq 8)
/\mbox{Span }([T_1],[T_2],[T_3]).
$$
Since $\Lambda$ is connected, the projections of its $(-2)$-spheres
must lie in the same lattice. This proves that for any component 
$\Lambda$, $L_\Lambda$ is contained in either
$$
\mbox{Span }([T_j],\xi_k|j=1,2,3, 1\leq k\leq 8)
/\mbox{Span }([T_1],[T_2],[T_3])
$$ 
or
$$
\mbox{Span }([T_j],\eta_k|j=1,2,3, 1\leq k\leq 8)
/\mbox{Span }([T_1],[T_2],[T_3]).
$$ 
Similarly, one can show that for any two distinct components 
$\Lambda,\Lambda^\prime$, the corresponding sublattices $L_\Lambda$, 
$L_{\Lambda^\prime}$ are orthogonal to each other.

Now let $\Lambda$ be a type (B) component, which consists of two
$(-2)$-spheres $C_1,C_2$ intersecting at a single point with 
tangency of order $2$. Because 
$$
({C}_1+{C}_2)^2=(-2)+2\cdot 2+(-2)=0,
$$ 
$\underline{C}_1+ \underline{C}_2$ must be $0$
and hence $L_\Lambda=\langle -2\rangle $ in this case.
Suppose $\Lambda$ is a type (C) component represented by a graph $\Gamma$
of type $\tilde{A}_n$, $\tilde{D}_n$, $\tilde{E}_6$, $\tilde{E}_7$
or $\tilde{E}_8$ listed in Figure 2, and let $\{C_i\}$ be the $(-2)$-spheres
corresponding to the vertices in $\Gamma$. Then there are weights $\{m_i\}$,
$m_i>0$, such that (1) $(\sum_i m_i C_i)^2=0$, (2) there exists a weight
$1=m_0\in\{m_i\}$ which has the property that if the corresponding vertex
(and the edge connecting to it) in $\Gamma$ is removed, the resulting graph 
is the Dynkin diagram for the root lattice
corresponding to $\Gamma$ (cf. \cite{BPV}, Lemma 2.12 (ii)). It follows 
easily that $L_\Lambda$ is isomorphic to the corresponding root lattice.

\end{proof}

With the preceding preparation, we give a proof of Theorem 1.8 next.

We assume $G\equiv\Z_p$ where $p=5$ or $7$. First, we observe that the 
main results in \cite{CK} (Theorem B, Theorem 3.1 and Prop. 3.7) 
concerning the fixed-point set structure of a symplectic $\Z_p$-action apply to the 
current situation, even though there was an additional assumption in \cite{CK}
that the symplectic $\Z_p$-action acts trivially on the second homology. This is
because the said additional assumption was mainly used to ensure that
the induced action of $G$ on each component of the union of $J$-holomorphic 
curves $\cup_i C_i$ leaves each $(-2)$-sphere in the component invariant  
if the component contains at least one fixed point, which is automatically
true in the current case. (Note that for $G\equiv\Z_5$ or $\Z_7$, the graphs
of type $\tilde{A}_n$, $\tilde{D}_n$, $\tilde{E}_6$, $\tilde{E}_7$
or $\tilde{E}_8$ listed in Figure 2 do not have any nontrivial $G$-symmetries
except for the case of $\tilde{A}_n$, in which $G$ acts freely so that the
corresponding component of $\cup_i C_i$ can not contain any fixed points 
of the $G$-action on $X_\alpha$.) With this understood, according to  
\cite{CK} the fixed points of $G$ can be divided into groups of the 
following types:

\begin{itemize}
\item [{(1)}] One fixed point with local representation 
$(z_1,z_2)\mapsto (\mu_p^k z_1,\mu_p^{-k}z_2)$ 
for some $k\neq 0\bmod p$, i.e., with local
representation contained in $SL_2(\C)$.
\item [{(2)}] Two fixed points with local representation 
$(z_1,z_2)\mapsto (\mu_p^{2k} z_1,\mu_p^{3k}z_2)$,
$(z_1,z_2)\mapsto (\mu_p^{-k} z_1,\mu_p^{6k}z_2)$ for some $k\neq 0
\bmod p$ respectively. (This type of fixed points occurs only when $p>5$.)
\item [{(3)}] Three fixed points, one with local representation 
$(z_1,z_2)\mapsto (\mu_p^k z_1,\mu_p^{2k}z_2)$ and the other two with 
local representation $(z_1,z_2)\mapsto (\mu_p^{-k} z_1,\mu_p^{4k}z_2)$ for
some $k\neq 0\bmod p$. 
\item [{(4)}] Four fixed points, one with local representation 
$(z_1,z_2)\mapsto (\mu_p^k z_1,\mu_p^k z_2)$ and the other three with 
local representation $(z_1,z_2)\mapsto (\mu_p^{-k} z_1,\mu_p^{3k}z_2)$ for
some $k\neq 0\bmod p$. 
\item [{($\Gamma$)}] The subset of fixed points which are contained in a
component $\Lambda$ of $\cup_i C_i$, where $\Lambda$ is of type (C) and
is represented by graph $\Gamma$ of type $\tilde{A}_n$, $\tilde{D}_n$, 
$\tilde{E}_6$, $\tilde{E}_7$ or $\tilde{E}_8$, such that at least one of
the $(-2)$-spheres in $\Lambda$ is fixed under the action. Note that
according to \cite{CK}, $n=-1\pmod{p}$ if $\Gamma$ is of type $\tilde{A}_n$,
and $n=4\pmod{p}$ if $\Gamma$ is of type $\tilde{D}_n$.
\item [{($T^2$)}] An embedded torus of self-intersection $0$. 
\end{itemize}

We shall consider the cases $p=5$ and $p=7$ separately.

\vspace{2mm}

{\it Case(a)} $p=5$. In this case, there are no type (2) fixed points. 
Moreover, by the assumption that both $\Theta_1,\Theta_2$ are nontrivial, 
and by Lemma 4.5 and Lemma 6.2, there are no type ($\Gamma$) fixed points 
unless $\Gamma$ is of type $\tilde{A}_4$ or $\tilde{D}_4$. The next lemma
further eliminates type ($\Gamma$) fixed points where $\Gamma$ is of type 
$\tilde{D}_4$.

\begin{lemma}
Let $G\subset \text{Aut}(E_8)$ be a subgroup of order $5$. There are no
sublattices of $E_8$ fixed under $G$, which are isomorphic either to a 
$D_4$-root lattice or to a direct sum of two copies of a $A_2$-root lattice.
\end{lemma}

\begin{proof}
By Lemma 4.5, there are two different integral $G$-representations
associated to the subgroup $G\subset\text{Aut}(E_8)$: $\Z[\Z_5]\oplus\Z^3$
and $\Z[\mu_5]^2$. The latter does not fix any vector in the lattice, so
we only need to consider the case of $\Z[\Z_5]\oplus\Z^3$.

By Carter \cite{Car} (Table 3, page 23), an element of order $5$ in 
$\text{Aut}(E_8)$ is uniquely determined up to conjugacy by the
characteristic polynomial. Hence a subgroup $G$ of order $5$ whose
corresponding integral $G$-representation is $\Z[\Z_5]\oplus\Z^3$
must be conjugate to the subgroup generated by the permutation
$$
e_1\mapsto e_2, e_2\mapsto e_3, e_3\mapsto e_4, e_4\mapsto e_5,
e_5\mapsto e_1, e_l\mapsto e_l, l=6,7,8,
$$  
where $\{e_1,\cdots,e_8\}$ is a standard basis of $\R^8$.
Clearly, the roots of $E_8$ which are fixed under the permutation 
can be put in two groups
$$
\Omega_1=\{\pm e_i\pm e_j|i\neq j, i,j=6,7,8\}
$$
and 
$$
\Omega_2=\{\frac{1}{2}\sum_{i=1}^8\epsilon_ie_i|\epsilon_1=\cdots=
\epsilon_5, \prod_{i=1}^8\epsilon_i=1\}.
$$
Note that for any roots $r_1,r_2\in\Omega_1$, $(r_1,r_2)=0$ if and only
if $r_1=\pm (e_i+e_j)$ and $r_2=\pm (e_i-e_j)$ (or vice versa), and 
$(r_1,r_2)\neq 0$ for any $r_1,r_2\in\Omega_2$.

With these preparations, we shall prove next that there are no sublattices
of $E_8$ isomorphic to a $D_4$-root lattice that are fixed under $G$. 
To see this, note that amongst the three roots represented by the vertices
other than the central one in a $D_4$-Dynkin diagram, exactly two of them must
belong to $\Omega_1$, which are of the form $\pm (e_i+e_j)$,
$\pm (e_i-e_j)$ for some $i\neq j$, $i,j=6,7,8$. On the other hand,
a root $r=\frac{1}{2}\sum_{i=1}^8\epsilon_i e_i\in\Omega_2$ is orthogonal
to $\pm (e_i+e_j)$ if and only if $\epsilon_i=-\epsilon_j$. But such a root
certainly is not orthogonal to $\pm (e_i-e_j)$. Our claim follows easily.

It remains to show that $G$ can not fix a direct sum of two copies of a
$A_2$-root lattice. To see this, let $r_1,r_2$ be the two roots generating
the first copy, and let $r_3,r_4$ generate the second copy. Then note first
that one of the $r_i$'s must belong to $\Omega_2$. Assume it is $r_1$ without
loss of generality. Then $r_3,r_4$, both being orthogonal to $r_1$, must 
belong to $\Omega_1$. Without loss of generality we may only consider the case 
$r_3=e_6-e_7$ and $r_4=e_7-e_8$. The root $r_1$, being orthogonal to both
$r_3,r_4$, must be $\pm \frac{1}{2}(e_1+\cdots +e_8)$. But then the root
$r_2$, which has the property that $(r_1,r_2)=-1$, can not be possibly
orthogonal to both $r_3,r_4$. A contradiction. The other cases are analogous,
and this finishes the proof of the lemma.

\end{proof}

We remark that there are sublattices of $E_8$ isomorphic to a $A_4$-root
lattice which are fixed under $G$. For example, the following $4$ roots
$$
e_6-e_7, -\frac{1}{2}(e_1+\cdots + e_5+e_6-e_7-e_8),
\frac{1}{2}(e_1+\cdots +e_5-e_6-e_7+e_8), e_6+e_7
$$
generate a $A_4$-root lattice which is fixed under $G$.

With the preceding understood, let $u,v,w$ and $A$ be the number of groups
of fixed points of $G$ of type (1), (3), (4) and ($\tilde{A}_4$)
respectively. We will next determine the possibilities of $u,v,w$ and $A$
using the Lefschetz fixed point theorem and the $G$-signature theorem
(i.e. Theorem 3.6). Note that since a fixed torus of self-intersection
$0$ makes no contribution in the calculation with the Lefschetz fixed point 
theorem and the $G$-signature theorem, we will ignore it in the consideration.
The number of such components in the fixed-point set will be determined later
by the number of cyclotomic summands in the integral representation on the
middle homology. 

To this end, recall that the total signature defect of a group of fixed points
of type (1), (3), (4) is $4$, $-8$, and $-4$ respectively (cf. Lemma 3.8
of \cite{CK}). For the total signature defect of a group of fixed points of
type ($\tilde{A}_4$), we note that such a component $\Lambda$ of $\cup_i C_i$
contains exactly $1$ fixed $(-2)$-sphere plus $3$ isolated fixed points of
local representation $(z_1,z_2)\mapsto (\mu_p^{k} z_1,\mu_p^{kq}z_2)$ for
some $k\neq 0\pmod{p}$ (where $p=5$), with $q=1,2,3$ respectively (cf.
\cite{CK}). The signature defect for each of the isolated fixed points
is correspondingly given by
$$
I_{p,q}=\sum_{k=1}^{p-1}\frac{(1+\mu_p^k)(1+\mu_p^{kq})}
{(1-\mu_p^k)(1-\mu_p^{kq})}.
$$
(Observe the relation $I_{p,q}=-I_{p,-q}$.) It follows easily that
$I_{5,1}=-4$ and $I_{5,2}=I_{5,3}=0$ (cf. \cite{CK}, Appendix). Hence
the total signature defect of a group of fixed points of type ($\tilde{A}_4$)
is 
$$
-4+0+0+\frac{5^2-1}{3}\cdot (-2)=-20.
$$

For $i=1,2$, let $\Z[\Z_5]^{r_i}\oplus\Z^{t_i}\oplus\Z[\mu_5]^{s_i}$
be the integral $G$-representation associated to $\Theta_i:G\rightarrow
\text{Aut}(E_8)$. Then by Lemma 4.5, there are the following three 
possibilities if we assume both $\Theta_1,\Theta_2$ are nontrivial:
$(r_1,t_1,s_1)=(r_2,t_2,s_2)=(1,3,0)$, $(r_1,t_1,s_1)=(r_2,t_2,s_2)
=(0,0,2)$ and  $(r_1,t_1,s_1)=(1,3,0)$, $(r_2,t_2,s_2)=(0,0,2)$. 
The $G$-signature theorem as stated in Theorem 3.6 
and the Lefschetz fixed point theorem (cf. Theorem 3.4)
give rise to the following equations:
$$
\left\{\begin{array}{cc}
1+3\cdot 2+t_1-s_1+t_2-s_2+1 = & u+3v+4w+5A\\
p\cdot (-r_1-t_1-r_2-t_2) = & -16+ 4u-8v-4w-20A\;
(\mbox { with } p=5).\\
\end{array} \right .
$$
The solution to the above system of equations is
$$
(u,v)=\left\{\begin{array}{cc}
(2-w+A,4-w-2A) & \mbox{ if } (r_1,t_1,s_1)=(r_2,t_2,s_2)=(1,3,0)\\
(3-w+A,2-w-2A) & \mbox{ if } (r_1,t_1,s_1)=(1,3,0), (r_2,t_2,s_2)=(0,0,2)\\
\end{array} \right .
$$
and $(u,v,w,A)=(4,0,0,0)$ if $(r_1,t_1,s_1)=(r_2,t_2,s_2)=(0,0,2)$. 

We shall further analyze the fixed-point set with help of the $G$-signature 
theorem as stated in Theorem 3.5 and with help of the $G$-index theorem
for Dirac operators as stated in Lemma 3.8. 

We consider first the cases where $A=0$, i.e., there are no type 
($\tilde{A}_4$) fixed points. To apply the $G$-signature theorem 
in Theorem 3.5, we fix a $g\in G$ and recall, with $p=5$ below, that 
each isolated fixed point $m$ of $g$ is associated with a pair of 
integers $(a_m,b_m)$, where $0<a_m,b_m<p$, such that the action of $g$ 
on the tangent space at $m$ is given by the complex linear transformation 
$(z_1,z_2)\mapsto (\mu_p^{a_m}z_1,\mu_p^{b_m}z_2)$, and 
moreover, the contribution to $\text{Sign }(g,X_\alpha)$ from $m$ is given by
$$
\delta_m=-\cot(\frac{a_m\pi}{p})\cdot\cot(\frac{b_m\pi}{p}).
$$
Now divide the fixed points of $g$ into three groups I, II, III according 
to their local representations: group I consists of fixed points with 
local representation $(z_1,z_2)\mapsto (\mu_p^k z_1,\mu_p^{k} z_2)$ for 
some $k\neq 0\pmod{p}$, group II consists of fixed points with 
$(z_1,z_2)\mapsto (\mu_p^k z_1,\mu_p^{2k} z_2)$ for some $k\neq 0\pmod{p}$, 
and group III consists of fixed points with 
$(z_1,z_2)\mapsto (\mu_p^k z_1,\mu_p^{4k} z_2)$ for some $k\neq 0\pmod{p}$. 
Then one observes that $\delta_m$ has only two possible 
values for the fixed points in each of the groups I, II, III. For group I, 
the values are $-\cot^2(\frac{\pi}{5})$, $-\cot^2(\frac{2\pi}{5})$, for 
group II, the values are $-\cot(\frac{\pi}{5})\cot(\frac{2\pi}{5})$,
$\cot(\frac{\pi}{5})\cot(\frac{2\pi}{5})$, and for group III, the 
values are $\cot^2(\frac{\pi}{5})$, $\cot^2(\frac{2\pi}{5})$. We let
$x_1,x_2$, $y_1,y_2$ and $z_1, z_2$ be the number of fixed points at which
$\delta_m$ takes these values respectively. 

By the $G$-signature theorem as stated in Theorem 3.5, we have
\begin{eqnarray*}
s_1+s_2-t_1-t_2=\text{Sign }(g,X_\alpha) & = & -x_1\cot^2(\frac{\pi}{5})
                                 -x_2\cot^2(\frac{2\pi}{5})\\
                           &   & -y_1\cot(\frac{\pi}{5})\cot(\frac{2\pi}{5})
                                 +y_2\cot(\frac{\pi}{5})\cot(\frac{2\pi}{5})\\
                           &   &   +z_1\cot^2(\frac{\pi}{5})
                                 +z_2\cot^2(\frac{2\pi}{5}).\\
\end{eqnarray*}
On the other hand, if we replace $g$ by $g^2$, $\delta_m$ will 
correspondingly be switched between the two values it assumes, and 
consequently, we have
\begin{eqnarray*}
s_1+s_2-t_1-t_2=\text{Sign }(g^2,X_\alpha) & = & -x_1\cot^2(\frac{2\pi}{5})
                                 -x_2\cot^2(\frac{\pi}{5})\\
                           &   & y_1\cot(\frac{\pi}{5})\cot(\frac{2\pi}{5})
                                 -y_2\cot(\frac{\pi}{5})\cot(\frac{2\pi}{5})\\
                           &   &   +z_1\cot^2(\frac{2\pi}{5})
                                 +z_2\cot^2(\frac{\pi}{5}).\\
\end{eqnarray*}
Combining these two equations, one obtains
\begin{eqnarray*}
[(z_1-z_2)-(x_1-x_2)]\cot^2(\frac{\pi}{5})-2(y_1-y_2)
\cot(\frac{\pi}{5})\cot(\frac{2\pi}{5})\\
-[(z_1-z_2)-(x_1-x_2)]\cot^2(\frac{2\pi}{5})=0.
\end{eqnarray*}

\begin{lemma}
$\cot(\frac{\pi}{5})/\cot(\frac{2\pi}{5})$ satisfies the 
algebraic equation $t^2-4t-1=0$, which is irreducible over $\Q$.
\end{lemma}

\begin{proof}
We start with the equation $1+\mu_5+\cdots +\mu_5^4=0$, from which one sees 
that $\cos(\frac{\pi}{5})$ satisfies $4t^2-2t-1=0$, and hence 
$\cos(\frac{\pi}{5})=(1+\sqrt{5})/4$.

Now observe that 
$$
\frac{\cot(\frac{2\pi}{5})}{\cot(\frac{\pi}{5})}=1-\frac{1}
{2\cos^2(\frac{\pi}{5})}.
$$
Using the fact that $\cos(\frac{\pi}{5})=(1+\sqrt{5})/4$, one can check 
that $\cot(\frac{\pi}{5})/\cot(\frac{2\pi}{5})$ is a solution of $t^2-4t-1=0$,
which is clearly irreducible over $\Q$.

\end{proof}
 
The preceding lemma implies the following relations 
$$
(z_1-z_2)-(x_1-x_2)=c,\; y_1-y_2=2c \mbox{ for some }c\in\Z.
$$
Now observe that type (1) fixed points contribute exclusively to $z_1$
or $z_2$, and we have $u=z_1+z_2$, and on the other hand, a group of type
(3) or type (4) fixed points contributes nontrivially to $x_1$ (resp. $x_2$)
if and only if it contributes nontrivially to $y_1$ (resp. $y_2$), and
we have 
$$
x_1=2v_1+w_1, \;\;x_2=2v_2+w_2, \;\;y_1=v_1+3w_1, \;\;y_2=v_2+3w_2
$$
where $v=v_1+v_2$ and $w=w_1+w_2$. From these equations we obtain
$$
2(z_1-z_2)=2(x_1-x_2)+(y_1-y_2)=5(v_1+w_1-v_2-w_2).
$$
Since $|z_1-z_2|\leq z_1+z_2=u<5$ in all the cases, we must have 
$$
z_1-z_2=v_1+w_1-v_2-w_2=0. 
$$
In particular, note that $u=z_1+z_2=2z_1$ is an even number.

The solutions which satisfy the above constraints are given below
(up to changing from $g$ to $g^2$)
\begin{itemize}
\item [{(a)}] $x_1=x_2=4, y_1=y_2=2, z_1=z_2=1, \mbox{ and } (u,v,w)=(2,4,0)$,
\item [{(b)}] $x_1=x_2=3, y_1=y_2=4, z_1=z_2=0, \mbox{ and } (u,v,w)=(0,2,2)$,
\item [{(c)}] $x_1=4$, $x_2=2$, $y_1=2$, $y_2=6$, $z_1=z_2=0$, and 
$(u,v,w)=(0,2,2)$,
\item [{(d)}] $x_1=2$, $x_2=1$, $y_1=1$, $y_2=3$, $z_1=z_2=1$, and 
$(u,v,w)=(2,1,1)$.
\end{itemize}

We next use Lemma 3.8 and Theorem 3.9 to rule out the cases (a), (b)  where 
$$
x_1-x_2=y_1-y_2=z_1-z_2=0 \mbox{ and } x_1-z_1=3.
$$

Observe that by the formula for the ``Spin-number'' in Lemma 3.8, 
the contribution to $\text{Spin }(g,X_\alpha)$ from a fixed point $m$ is 
$$
\nu_m=-(-1)^{k(g,m)}\cdot\frac{1}{4}\cdot\csc(\frac{a_m\pi}{5})
\csc(\frac{b_m\pi}{5}),
$$
where $0<a_m,b_m<5$ and $k(g,m)\cdot 5=2r_m+a_m+b_m$ for some $0\leq r_m<5$.
One can check that $\nu_m$ takes values $-\frac{1}{4}\csc^2(\frac{\pi}{5})$, 
$-\frac{1}{4}\csc^2(\frac{2\pi}{5})$ if $m$ belongs to group I; for group
II, the values of $\nu_m$ are 
$\frac{1}{4}\csc(\frac{\pi}{5})\csc(\frac{2\pi}{5})$,
$-\frac{1}{4}\csc(\frac{\pi}{5})\csc(\frac{2\pi}{5})$, and for group III, 
the values are $\frac{1}{4}\csc^2(\frac{\pi}{5})$, 
$\frac{1}{4}\csc^2(\frac{2\pi}{5})$.
The number of fixed points at which $\nu_m$ takes these values is 
$x_1, x_2, y_1,y_2,z_1,z_2$ respectively.

With the above understood, for the cases (a), (b) we obtain from Lemma 3.8
\begin{eqnarray*}
\text{Spin }(g,X_\alpha) & = & -\frac{x_1}{4}\csc^2(\frac{\pi}{5})
                        -\frac{x_2}{4}\csc^2(\frac{2\pi}{5})\\
                  &   & +\frac{y_1}{4}\csc(\frac{\pi}{5})\csc(\frac{2\pi}{5})
                        -\frac{y_2}{4}\csc(\frac{\pi}{5})\csc(\frac{2\pi}{5})\\
                  &   & +\frac{z_1}{4}\csc^2(\frac{\pi}{5})+
                         \frac{z_2}{4}\csc^2(\frac{2\pi}{5})\\
                  & = & \frac{z_1-x_1}{4}\csc^2(\frac{\pi}{5})+
                        \frac{z_2-x_2}{4}\csc^2(\frac{2\pi}{5})\\
                  & = & z_1-x_1\\
                  & = & -3\\
\end{eqnarray*}
because $z_1-x_1=z_2-x_2$ and 
$\csc^2(\frac{\pi}{5})+\csc^2(\frac{2\pi}{5})=4$. 

On the other hand, there are integers $d_0,\cdots, d_4$ such that
$$
\text{Spin }(g,X_\alpha)=d_0+d_1\mu_5+\cdots +d_4\mu_5^4.
$$
Since $1+t+\cdots +t^4=0$ is irreducible over $\Q$, one must have
$$
d_0+3=d_1=\cdots =d_4.
$$
With the fact that the index of the Dirac operator on $X_\alpha$, which is
given by the sum $d_0+\cdots +d_4$, equals $-\text{Sign}(X_\alpha)/8=2$, we
obtain $d_0=-2$ and $d_1=\cdots=d_4=1$. By Fang's theorem (cf. Theorem 3.9),
the Seiberg-Witten invariant 
$$
SW_{X_\alpha}(0)=0 \pmod {5}
$$ 
for the trivial $Spin^\C$-structure on $X_\alpha$. (Note that the trivial 
$Spin^\C$-structure on $X_\alpha$ is a $G$-$Spin^\C$ structure because 
by Lemma 3.8, the action of $G$ is spin.) However, this is a contradiction, 
because by construction $SW_{X_\alpha}(0)=1$, cf. Section 2. This proves 
our claim regarding the cases (a), (b). 

For case (c), a similar calculation shows that
$$
\text{Spin }(g,X_\alpha)=-2+2\mu_5^2+2\mu_5^3
$$
(i.e., $d_0=-2$, $d_2=d_3=2$ and $d_1=d_4=0$), which does not violate 
Fang's theorem (cf. Theorem 3.9). With $d_0=-2$, this set of fixed-point
data does not violate Theorem 3.10 either. 

For case (d), we have by a similar calculation that
$$
\text{Spin }(g,X_\alpha)=\mu_5^2+\mu_5^3, 
$$
which violates Fang's theorem, and so it is eliminated. 

To finish the analysis for the cases where $A=0$, it remains to check case (c)
against Theorem 3.11, a constraint coming from the Kirby-Siebenmann and the
Rochlin invariants. One finds easily from the fixed-point set structure that the
corresponding $4$-manifold with boundary $N$ has $14$ boundary components:
there are six $L(5,1)$, six $L(5,2)$, and two $L(5,3)$. By Corollary 2.24 in \cite{Sa},
the Rochlin invariant of $L(5,1)$,  $L(5,2)$, and $L(5,3)$ equals $4$, $0$, and $0$
respectively. On the other hand, one finds easily that $\text{Sign}(N)=-8$. The equation
in Theorem 3.11 becomes 
$$
8\cdot \text{ks}(N)\equiv -8+6\times 4+6\times 0+2\times 0 \pmod{16}, 
$$
which implies that $\text{ks}(N)=0$. Hence case (c) can not be ruled out by Theorem 3.11. 

Now we consider the cases where there are type ($\tilde{A}_4$) fixed points,
i.e., $A\neq 0$. Then by the analysis in \cite{CK} 
(cf. Lemma 3.6 of \cite{CK}), we see that there exists a $k\neq 0\pmod{p}$
such that the local representations of $g$ at the $3$ isolated points
in a group of type ($\tilde{A}_4$) fixed points are 
$$(z_1,z_2)\mapsto (\mu_p^{-3k} z_1,\mu_p^{-k} z_2),
(z_1,z_2)\mapsto (\mu_p^{3k} z_1,\mu_p^{3k} z_2),
(z_1,z_2)\mapsto (\mu_p^{-3k} z_1,\mu_p^{-k} z_2),
$$
and the local representation at the fixed $(-2)$-sphere is
$z\mapsto \mu_p^k z$ (here $p=5$). The key observation is that the total
contribution from a group of type ($\tilde{A}_4$) fixed points to
$\text{Sign }(g,X_\alpha)$ equals $-5$, which is independent of $g$. More
concretely, by a direct calculation the total contribution is
$$
-2\cot(\frac{-3k\pi}{5})\cdot\cot(\frac{-k\pi}{5})-
\cot^2(\frac{3k\pi}{5})-2\csc^2(\frac{k\pi}{5})=-5, \;\forall k.
$$
Consequently, one can similarly introduce the
numbers $x_1,x_2,y_1,y_2,z_1,z_2$ and $v_1,v_2$, $w_1,w_2$
for the fixed points of type (1), type (3),
or type (4) as in the case of $A=0$, and the same argument implies that
$$
z_1-z_2=v_1+w_1-v_2-w_2=0. 
$$
In particular, $u=z_1+z_2=2z_1$ is an even number. 

With Lemma 6.3, the solutions (up to changing from $g$ to $g^2$) 
which satisfy these constraints are
\begin{itemize}
\item [{(i)}] $x_1=x_2=y_1=y_2=0, z_1=z_2=2, \mbox{ and } (u,v,w,A)=(4,0,0,2)$,
\item [{(ii)}] $x_1=2$, $x_2=1$, $y_1=1$, $y_2=3$, $z_1=z_2=1$, and 
$(u,v,w,A)=(2,1,1,1)$.
\item [{(iii)}] $x_1=x_2=y_1=y_2=0, z_1=z_2=2, 
\mbox{ and } (u,v,w,A)=(4,0,0,1)$.
\end{itemize}

Next we use Lemma 3.8 and Fang's theorem (cf. Theorem 3.9) to examine these
fixed-point data. To this end, we need to determine the possible values 
of the total contribution of a group of type ($\tilde{A}_4$) fixed points
to the ``Spin-number'' $\text{Spin }(g,X_\alpha)$. A direct calculation
shows that for $k=1,4$, the total contribution is 
$$
-\frac{2}{4}\csc(\frac{2\pi}{5})\cdot \csc(\frac{4\pi}{5})
-\frac{1}{4}\csc^2(\frac{3\pi}{5})-\frac{-2}{4}\csc(\frac{\pi}{5})\cdot
\cot(\frac{\pi}{5})=0,
$$
and for $k=2,3$, it is 
$$
\frac{2}{4}\csc(\frac{4\pi}{5})\cdot \csc(\frac{3\pi}{5})
-\frac{1}{4}\csc^2(\frac{\pi}{5})+\frac{-2}{4}\csc(\frac{2\pi}{5})\cdot
\cot(\frac{2\pi}{5})=0.
$$
As an immediate consequence we obtain that the ``Spin-number''
$$
\text{Spin }(g,X_\alpha)=\left\{\begin{array}{cc}
2 & \mbox{ in case (i)}\\
\mu_5^2+\mu_5^3 & \mbox{ in case (ii)}\\
2  & \mbox{ in case (iii)}\\
\end{array} \right .
$$
Note that case (ii) violates Fang's theorem, hence is eliminated. 
However, the remaining cases (i) and (iii) can not be ruled out by
Theorem 3.10 (in both cases $d_0=2$). 

It remains to show that the number of fixed tori in the fixed-point set
is bounded from above by one half of the number of copies of $\Z[\mu_5]$
in the associated integral $G$-representation of $\Theta=(\Theta_1,
\Theta_2): G\rightarrow \text{Aut}(E_8\oplus E_8)$. To see this, let
$H_2(X_\alpha;\Z)=\Z[\Z_5]^r\oplus \Z^t\oplus \Z[\mu_5]^s$ be the 
decomposition into summands of regular, trivial and cyclotomic types.
Since the classes $[T_1]$, $[T_2]$, $[T_3]$ are fixed under the $G$-action
on $X_\alpha$, we see that $\Z[\mu_5]^s$ is orthogonal to  $[T_1]$, $[T_2]$, 
and $[T_3]$, hence is contained in 
$$
\mbox{Span }([T_j],\xi_k,\eta_k|j=1,2,3, 1\leq k\leq 8).
$$ 
(Here $\xi_k,\eta_k$ are the classes in $H_2(X_\alpha;\Z)$ which correspond 
to the two standard bases of the $-E_8$ form defined in Lemma 4.2.)
The number $s$ is bounded from above by the number of copies of $\Z[\mu_5]$
in the associated integral $G$-representation of $\Theta$ follows from 
Proposition 4.6 plus the fact that $\Z[\mu_5]^s$ is mapped injectively into 
$$
\mbox{Span }([T_j],\xi_k,\eta_k|j=1,2,3, 1\leq k\leq 8)/
\mbox{Span }([T_1],[T_2],[T_3])
$$
under the quotient map. Our claim then follows from Edmonds' result
(cf. Prop. 3.1). The case of Theorem 1.8 where $p=5$ follows easily. 

\vspace{3mm}

{\it Case(b)} $p=7$. In this case, there are no type ($\Gamma$) fixed 
points by Lemma 4.5 and Lemma 6.2, because by the assumption $\Theta_1,
\Theta_2$ are nontrivial. Moreover, there are no fixed tori of 
self-intersection $0$ either by Proposition 4.6. The next lemma eliminates
the possibility of having type (4) fixed points.

\begin{lemma}
Let $G\subset \text{Aut}(E_8)$ be a subgroup of order $7$. There are no
sublattices of $E_8$ fixed under $G$, which are isomorphic to a $A_2$-root 
lattice.
\end{lemma}

\begin{proof}
By Carter \cite{Car} (Table 3, page 23), an element of order $7$ in 
$\text{Aut}(E_8)$ is uniquely determined up to conjugacy by the
characteristic polynomial. Hence a subgroup $G$ of order $7$, with the
corresponding integral $G$-representation being $\Z[\Z_5]\oplus\Z^3$
(cf. Lemma 4.5), must be conjugate to the subgroup generated by the 
permutation
$$
e_1\mapsto e_2, e_2\mapsto e_3, \cdots, e_6\mapsto e_7, 
e_7\mapsto e_1, \mbox{ and } e_8\mapsto e_8,
$$  
where $\{e_1,\cdots,e_8\}$ is a standard basis of $\R^8$. The only roots
which are fixed under the permutation are
$$
r=\pm \frac{1}{2}(e_1+e_2+\cdots +e_7+e_8),
$$
which do not generate a $A_2$-root lattice. The lemma follows.

\end{proof}

With the preceding understood, we shall next determine the number of groups
of type (1), type (2), and type (3) fixed points, which is denoted by
$u,v,w$ respectively. By Lemma 3.8 of \cite{CK}, the total signature 
defect of each of such groups is $10$, $-8$, and $2$ respectively. The 
Lefschetz fixed point theorem and the $G$-signature theorem (as in 
Theorem 3.6) give rise to the following equations
$$
\left\{\begin{array}{cc}
1+3\cdot 2+1+1+1 = & u+2v+3w\\
p\cdot (-1-1-1-1) = & -16+10u-8v+2w\; (\mbox { with } p=7).\\
\end{array} \right .
$$
The solutions are $(u,v,w)=(0,2,2), (1,3,1), (2,4,0)$. 

We examine these data with the $G$-signature theorem as in Theorem 3.5
and the $G$-index theorem for Dirac operators in Lemma 3.8.  

Let $\delta_1,\delta_2,\delta_3$ be the total contributions 
to $\text{Sign}(g,X_\alpha)$ of a group
of fixed points of type (1), type (2), type (3) respectively. With a
direct calculation we list all the possible values of them 
(in approximations) below, taken at $k=1,2,3$ respectively.
\begin{itemize}
\item $\delta_1=4.31194, \;\;\;0.63596, \;\;\;0.05210$,
\item $\delta_2=-4.49396, \;\;\;-1.10992, \;\;\;1.60388$,
\item $\delta_3=-2.60388, \;\;\;3.49396, \;\;\;0.10992$.
\end{itemize}

The cases where $(u,v,w)=(1,3,1), (2,4,0)$ can be eliminated as follows.
Consider the case of $(u,v,w)=(1,3,1)$ first. There are $3$ groups of type
(2) fixed points. If all three values of $\delta_2$ are assumed, then 
because the sum of these three values of $\delta_2$ equals $-4$ and 
$\text{Sign}(g,X_\alpha)=-2$, the sum of the values of $\delta_1$ and
$\delta_3$ must equal $2$, which is easily seen impossible by examining 
the possible values of $\delta_1$ and $\delta_3$. If not all three values
of $\delta_2$ are assumed, then there must be one value of $\delta_2$ 
which is assumed by at least $2$ groups of type (2) fixed points. By 
changing $g\in G$ to a suitable power of $g$, we may assume this value of
$\delta_2$ is $-4.49396$. Then the sum of the rest of the values, one for
each of $\delta_1,\delta_2,\delta_3$, must be $2\times 4.49396-2=6.98792$.
But this is easily seen impossible by examining the possible values.
This ruled out the case where $(u,v,w)=(1,3,1)$. By a similar argument,
the case where $(u,v,w)=(2,4,0)$ can be also ruled out. 

It remains to examine the case where $(u,v,w)=(0,2,2)$. Observe first
that for each value of $\delta_2$, there is a unique value of $\delta_3$
such that the sum of the two values equals $-1$. It follows easily from
this that the fixed points of $g$ can be divided into two groups, where
each group consists of $5$ points with local representations 
$(z_1,z_2)\mapsto (\mu_p^{2k}z_1,\mu_p^{3k}z_2)$,
$(z_1,z_2)\mapsto (\mu_p^{-k}z_1,\mu_p^{-k}z_2)$,
$(z_1,z_2)\mapsto (\mu_p^{2k}z_1,\mu_p^{4k}z_2)$,
$(z_1,z_2)\mapsto (\mu_p^{-2k}z_1,\mu_p^kz_2)$,
$(z_1,z_2)\mapsto (\mu_p^{-2k}z_1,\mu_p^kz_2)$
respectively, for some $k\neq 0\pmod{p}$. The question is whether the
number $k\pmod{p}$ (which is only determined up to a sign) must be the
same for the two groups. We will show that it must be the same using
Lemma 3.8.

To this end, we let $\nu_2,\nu_3$ be the total contributions to
the ``Spin-number'' $\text{Spin}(g,X_\alpha)$ of a group
of fixed points of type (2), type (3) respectively. With a
direct calculation we list all the possible values of them 
(in approximations) below, taken at $k=1,2,3$ respectively.
\begin{itemize}
\item $\nu_2=-1, \;\;\;-1, \;\;\;-1,$
\item $\nu_3=-0.44504, \;\;\;-1.80194, \;\;\;1.24698$.
\end{itemize}
On the other hand, $\text{Spin}(g,X_\alpha)=\sum_{l=0}^{p-1}d_l\mu_p^l$
where $d_0$ is even and $d_l=d_{p-l}$ for $l\neq 0$. With the observation
that
$$
2\cos(\frac{2\pi}{7})=1.24698,\;\;\; 2\cos(\frac{4\pi}{7})=-0.44504,\;\;\;
2\cos(\frac{6\pi}{7})=-1.80194,
$$
we can easily conclude that in  
$\text{Spin}(g,X_\alpha)=\sum_{l=0}^{p-1}d_l\mu_p^l$, $d_0=-2$ and
$d_1,\cdots, d_6$ contain two $0$'s and four $1$'s if $\nu_3$ assumes
two distinct values, and $d_1,\cdots, d_6$ contain four $0$'s and two $2$'s 
if $\nu_3$ assumes only one value. The former case violates Fang's theorem
(cf. Theorem 3.9), so it can not occur. This proves that the number 
$k\pmod{p}$ in the local representations must be the same for the two 
different groups of fixed points. Finally, we note that since $d_0=-2$ in 
the ``Spin-number'' $\text{Spin}(g,X_\alpha)=\sum_{l=0}^{p-1}d_l\mu_p^l$,
the remaining case can not be ruled out by Theorem 3.10. Moreover, by a 
similar argument as in the case of $p=5$, one can check easily that Theorem 3.11 
is not violated either. The proof for the case of $p=7$ in Theorem 1.8 is completed.

\vspace{0.5cm}

\noindent{\bf Acknowledgment}

\vspace{3mm}

We are grateful to Jin Hong Kim and Nobuhiro Nakamura for pointing out some 
errors in an earlier version of this article. We are also indebted to
Stefan Friedl for stimulating discussions (with the first author) which 
have triggered a substantial revision of this work. Finally, we wish to thank an 
anonymous referee whose extensive comments have led to a much improved 
exposition.

\vspace{1cm}

{\Small 
W. Chen: Department of Math. and Stat., University of Massachusetts,
Amherst, MA 01003.\\
{\it e-mail:} wchen@math.umass.edu

S. Kwasik: Mathematics Department, Tulane University,
New Orleans, LA 70118. \\
{\it e-mail:} kwasik@math.tulane.edu\\
}

\end{document}